\documentclass[a4paper,8pt]{article}
\usepackage{latexcad}
\usepackage{amsmath}
\usepackage{graphicx}
\usepackage{amsfonts}
\usepackage{amssymb}
\newtheorem{thm}{Theorem}[section]
\newtheorem{prop}[thm]{Proposition}
\newtheorem{cor}[thm]{Corollary}
\newtheorem{lem}[thm]{Lemma}
\newtheorem{dfn}[thm]{Definition}

\begin{document}

\title{The compositional construction of Markov processes II}
\author{L. de Francesco Albasini
\and N. Sabadini
\and R.F.C. Walters}
\maketitle

%
%
%
%
%
%
%
%
%
%
%
%
%
%
%
%
%
%
%
%
%
%
%
%
%
%
%

\section{Introduction}

In \cite{dFSW09a} we introduced a notion of Markov automaton, together with
parallel operations which permit the compositional description of Markov
processes. We illustrated by showing how to describe a system of $n$ dining
philosophers (with $12^{n}$ states), and we observed that Perron-Frobenius
theory yields a proof that the probability of reaching deadlock tends to one
as the number of steps goes to infinity. In this paper we add sequential
operations to the algebra (and the necessary structure to support them)
following analogous developments in \cite{KSW00, dFSW09b}. The extra
operations permit the description of hierarchical systems, and ones with
evolving geometry. We illustrate our algebra by describing a system called
Sofia's Birthday Party, originally introduced in \cite{KSW00}.

There is a huge literature on probabilistic and weighted automata,
transducers, and process calculi (see for example,
\cite{BEP97,H96,LSV03,R63,SV2004,DKV09}). However the model of \cite{dFSW09a}
is distinguished from the others in the following ways:

\begin{enumerate}
\item [(i)]In other probabilistic automata models (\cite{R63},\cite{BGC09})
the sum of probabilities of actions out of a given state \emph{with a given
label} is 1. This means that the probabilities are conditional on the
existence of the label. In our model the sum of probabilities of actions out
of a given state \emph{for all labels} \ is 1. The reason for this is further
explained by the next point.

\item[(ii)] Our weighted automata are very close mathematically to weighted
transducers but conceptually very different. Instead of modelling devices
which translate input to output, the idea is we model devices with number of
parallel interfaces, and when a transition occurs in the device this induces a
transition on each of the interfaces (the interfaces are part of the device).
In order to have binary operations of composition the interfaces are divided
into left and right interfaces. The notions of initial and final states are
also generalized in our notion of weighted automaton to become instead
functions into the state space. These sequential interfaces are not to be
thought of as initial and final states, but hooks into the state space at
which a behaviour may enter or leave the device. The application of our
weighted automata is to concurrent hierarchical and distributed systems rather
than language recognition or processing. In \cite{dFSW09b} we show how data types and
also state on the connectors (shared variables) may be added to our model.

\item[(iii)] For many compositional models the communication is based on
process algebras like CCS, CSP. Instead ours models truly concurrent systems
with explicit network topologies. One of the key aspects is our algebra of
connectors: parallel and sequential ``wires", which give an hierarchical network
topology to expressions in the algebra.

\item[(iv)] Our automata both with respect to the sequential and with respect
to the parallel operations form the arrows of a symmetric monoidal categories,
with other well-known categorical properties (see, for example, \cite{CW87}).
The operations are in fact based on the operations available in the categories
of spans and cospans. The categorical operations have an associated geometry.
Important equations satisfied are the Frobenius and separable axioms
\cite{RSW05}.
\end{enumerate}

As we have said the novelty of this paper is the introduction of
\emph{sequential operations} to the algebra of \cite{dFSW09a}. This requires
extra structure to be added to the automata, namely sequential interfaces.

Further, in \cite{dFSW09b} weighted automata (where the weighting of a
transition is a non-negative real number) played a subsidiary role. However
for sequential composition weighted automata are more fundamental, since in
identifying states of two different automata it is the relative weight given
to decisions which is important rather than the probabilities. Technically
this appears in the fact that normalization is not compositional with respect
to sequential operations, whereas for parallel operations it is. For a summary
of recent work on the various kinds of weighted automata see \cite{DKV09}, and
in particular \cite{MM09}.

Another aim of the paper is to illustrate how hierarchical and mobile systems may
be modelled in this algebra, using the combined sequential and parallel 
operations. 
Given a set of automata $S$, let us denote the set of automata given as
expressions in terms of parallel operations in the automata $S$ as $\Pi(S)$
and given as expressions in terms of sequential operations as $\Sigma(S)$. Let
$E$ be the set of elementary automata with only one transition. Then any
automaton has a representation in $\Sigma(E)$. The dining philosopher problem
of \cite{dFSW09b} is described as an element of $\Pi\Sigma(E)$, that is of
communicating sequential systems. An element of $\Sigma\Pi\Sigma(E)$ is one in
which the parallel structure may evolve, and so on. The system Sofia's
Birthday Party is in $\Pi\Sigma\Sigma(E)$ but illustrates also the form of
systems of type $\Pi\Sigma\Pi\Sigma(E)$. The power of the sequential and
parallel operations is that they may be alternated, as the alternating
quantifiers in logic, or the alternation in alternating Turing machines, or
the alternating sums and products in State Charts.

\section{Weighted automata with parallel and sequential interfaces}

In this section we define weighted and Markov automata \emph{with sequential
and parallel interfaces}, which however we shall call just weighted and Markov
automata. The reader should be aware that the definitions of \cite{dFSW09b}
differ in lacking sequential interfaces. We also do not require here for
weighted automata the special symbols $\varepsilon$ and the condition that the
rows of the total matrix are strictly positive: we reserve those conditions
for what we now call \emph{positive weighted automata}.

Notice that in order to conserve symbols in the following definitions we shall
use the same symbol for the automaton, its state space and its family of
matrices of transitions, distinguishing the separate parts only by the font.

\begin{dfn}
Consider two finite alphabets $A$ and $B$, and two finite sets $X$ and $Y$.
A\emph{\ weighted automaton }$\mathbf{Q}$ \emph{with left parallel interface}
$A$, \emph{right parallel interface} $B$, \emph{top sequential interface} $X$,
and \emph{bottom sequential interface} $Y$, consists of a finite set $Q$ of
states, and an $A\times B$ indexed family $\mathsf{Q=(Q}_{a,b}\mathsf{)}%
_{(a\in A,b\in B)}$ of $Q\times Q$ matrices with non-negative real
coefficients, and two functions, $\gamma_{0}:X\rightarrow Q$, and $\gamma
_{1}:Y\rightarrow Q$. We denote the elements of the matrix $\mathsf{Q}_{a,b}$
by $[\mathsf{Q}_{a,b}]_{q,q^{\prime}}$ $(q,q^{\prime}\in Q)$.
\end{dfn}

We call the matrix $\mathcal{Q}=\sum_{a\in A,b\in B}\mathsf{Q}_{a,b}$ the
\emph{total matrix} of the automaton.

We will use a brief notation for the automata $\mathbf{Q}$ indicating its
interfaces, namely $\mathbf{Q}_{Y;A,B}^{X}$. We shall use the same symbols
$\gamma_{0}$, $\gamma_{1}$ for the sequential interface functions of any
automata, and we will sometimes refer to these functions as the sequential
interfaces. Notice that the terms `left', `right', `top' and `bottom' for the
interfaces have no particular semantic significance - they are chosen to be
semantically neutral in order not to suggest input, output, initial or final.

\begin{dfn}
A weighted automaton $\mathbf{Q}$ is \emph{positive} if the parallel
interfaces $A$ and $B$ contain special elements, the symbols $\varepsilon_{A}$
and $\varepsilon_{B}$, and satisfies the property that the row sums of the
matrix $\mathsf{Q}_{\varepsilon_{A},\varepsilon_{B}}$ are strictly positive.
\end{dfn}

For a positive weighted automaton the total matrix has strictly positive row sums.

\begin{dfn}
A \emph{Markov automaton} $\mathbf{Q}$\emph{\ with left interface} $A$,
\emph{right interface} $B$, \emph{top sequential interface} $X$, and
\emph{bottom sequential interface} $Y$, written briefly $\mathbf{Q}%
_{Y;A,B}^{X}$, is a positive weighted automaton satisfying the extra condition
that the row sums of the total matrix $\mathcal{Q}$ are all $1$. That is, for
all $q$
\[
\sum_{q^{\prime}}\sum_{a\in A,b\in B}\left[  \mathsf{Q}_{a,b}\right]
_{q,q^{\prime}}=1.
\]
\end{dfn}

For a Markov automaton we call $[\mathsf{Q}_{a,b}]_{q,q^{\prime}}$ the
\emph{probability of the transition from }$q$\emph{\ to }$q^{\prime}%
$\emph{\ with left signal }$a$ \emph{and right signal} $b$.

The idea is that in a given state various transitions to other states are
possible and occur with various probabilities, the sum of these probabilities
being $1$. The transitions that occur have effects, which we may think of a
\emph{signals}, on the two interfaces of the automaton, which signals are
represented by letters in the alphabets. We repeat that it is fundamental
\emph{not} to regard the letters in $A$ and $B$ as inputs or outputs, but
rather signals induced by transitions of the automaton on the interfaces. For
examples see section 2.

When both $A$ and $B$ are one element sets and $X=Y=\emptyset$ a Markov
automaton is just a Markov matrix.

\begin{dfn}
Consider a weighted automaton $\mathbf{Q}$ with interfaces $A$ and $B$. A
\emph{behaviour} of length $k$ of $\mathbf{Q}$ consists of two words of length
$k$, one $u=a_{1}a_{2}\cdots a_{k}$ in $A^{\ast}$ and the other $v=b_{1}%
b_{2}\cdots b_{k}$ in $B^{\ast}$ and a sequence of non-negative row vectors
\[
x_{0},\;x_{1}=x_{0}\mathsf{Q}_{a_{1},b_{1}},\;x_{2}=x_{1}\mathsf{Q}%
_{a_{2},b_{2}},\;\cdots,\;x_{k}=x_{k-1}\mathsf{Q}_{a_{k},b_{k}}.
\]
Notice that, even for Markov automata, $x_{i}$ is \emph{not} generally a
distribution of states; for example often $x_{i}=0.$
\end{dfn}

\begin{dfn}
The \emph{normalization} of a positive weighted automaton $\mathbf{Q}$,
denoted $\mathbf{N(Q)}$ is the Markov automaton with the same interfaces and
states, but with
\[
\left[  \mathsf{N(Q)}_{a,b}\right]  _{q,q^{\prime}}=\frac{\left[
\mathsf{Q}_{a,b}\right]  _{q,q^{\prime}}}{\sum_{q^{\prime}\in Q}\left[
\mathcal{Q}\right]  _{q,q^{\prime}}}=\frac{\left[  \mathsf{Q}_{a,b}\right]
_{q,q^{\prime}}}{\sum_{q^{\prime}\in Q}\sum_{a\in A,b\in B}\left[
\mathsf{Q}_{a,b}\right]  _{q,q^{\prime}}}.
\]
\end{dfn}

It is obvious that a weighted automaton $\mathbf{Q}$ is Markov if and only if
$\mathbf{Q=N(Q)}$.

\begin{dfn}
If\ $\mathbf{Q}$ is a weighted automaton and $k$ is a natural number, then the
\emph{automaton of }$k$\emph{\ step paths }in $\mathbf{Q}$, which we denote as
$\mathbf{Q}^{k}$ is defined as follows: the states of $\mathbf{Q}^{k}$ are
those of $\mathbf{Q}$; the sequential interfaces are the same; the left and
right interfaces are $A^{k}$ and $B^{k}$ respectively.\ If $u=(a_{1}%
,a_{2},\cdots,a_{k})\in A^{k}$ and $v=(b_{1},b_{2},\cdots,b_{k})\in B^{k}$
then
\[
\mathbf{(}\mathsf{Q}^{k})_{u,v}=\mathsf{Q}_{a_{1},b_{1}}\mathsf{Q}%
_{a_{2},b_{2}}\cdots\mathsf{Q}_{a_{k},b_{k}}.
\]
The definition for positive weighted automata requires in addition that
$\varepsilon_{A^{k}}=(\varepsilon_{A},\cdots,\varepsilon_{A}),\varepsilon
_{B^{k}}=(\varepsilon_{B},\cdots,\varepsilon_{B})$.
\end{dfn}

If $\mathbf{Q}$ is a weighted automaton and $u=(a_{1},a_{2},\cdots,a_{k})\in
A^{k}$, $v=(b_{1},b_{2},\cdots,b_{k})\in B^{k}$, then $\mathbf{[(}%
\mathsf{Q}^{k})_{u,v}]_{q,q^{\prime}}$ is the sum over all paths from $q$ to
$q^{\prime}$ with left signal sequence $u$ and right signal sequence $v$ of
the weights of paths, where the weight of a path is the product of the weights
of the steps.

\subsection{Graphical representation of weighted automata}

Although the definitions above are mathematically straightforward, in practice
a graphical notation is more intuitive. We may compress the description of an
automaton with parallel interfaces $A$ and $B$, which requires $A\times B$
matrices, into a single labelled graph, like the ones introduced in
\cite{KSW97b}. We indicate by describing some examples.

\subsubsection{An example}

Consider the automaton with parallel interfaces $\{a\}$, $\{b_{1}%
,b_{2}\}\times\{c\}$, sequential interfaces $\{x\}$, $\{y,z\}$; with states
$\{1,2,3\}$ sequential interface functions $x\mapsto1$ and $y,z\mapsto3$; and
transition matrices
\[
\mathsf{Q}_{a,(b_{1},c)}=\left[
\begin{array}
[c]{ccc}%
0 & 2 & 0\\
0 & 3 & 0\\
0 & 0 & 0
\end{array}
\right]  ,\;\mathsf{Q}_{a,(b_{2},c)}=\left[
\begin{array}
[c]{ccc}%
0 & 0 & 0\\
0 & 0 & 1\\
0 & 0 & 0
\end{array}
\right]  .
\]
This information may be put in the diagram:

%
%
%
%
%
%
%
%
%
%
%
%
%
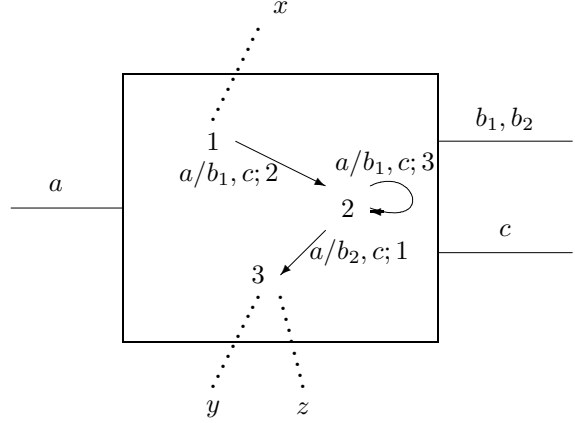
\centerline{\tt\setlength{\unitlength}{4.2pt}
\begin{picture}(58,48)
\thinlines
\drawframebox{28.0}{24.0}{28.0}{24.0}{}
\drawpath{14.0}{24.0}{4.0}{24.0}
\drawpath{42.0}{30.0}{54.0}{30.0}
\drawpath{42.0}{20.0}{54.0}{20.0}
\drawcenteredtext{22.0}{30.0}{$1$}
\drawcenteredtext{34.0}{24.0}{$2$}
\drawcenteredtext{26.0}{18.0}{$3$}
\drawcenteredtext{28.0}{42.0}{$x$}
\drawcenteredtext{22.0}{6.0}{$y$}
\drawcenteredtext{30.0}{6.0}{$z$}
\drawcenteredtext{8.0}{26.0}{$a$}
\drawcenteredtext{48.0}{32.0}{$b_1,b_2$}
\drawcenteredtext{48.0}{22.0}{$c$}
\drawvector{24.0}{30.0}{8.0}{2}{-1}
\drawvector{32.0}{22.0}{4.0}{-1}{-1}
\path
(36.0,26.0)(36.0,26.0)(36.11,26.04)(36.23,26.11)(36.34,26.15)(36.47,26.2)(36.58,26.24)(36.69,26.27)(36.8,26.29)(36.91,26.33)
\path
(36.91,26.33)(37.01,26.34)(37.12,26.36)(37.23,26.38)(37.33,26.4)(37.44,26.4)(37.54,26.4)(37.63,26.4)(37.73,26.4)(37.83,26.4)
\path
(37.83,26.4)(37.93,26.38)(38.01,26.38)(38.11,26.36)(38.19,26.34)(38.27,26.31)(38.36,26.29)(38.44,26.27)(38.52,26.25)(38.59,26.2)
\path
(38.59,26.2)(38.68,26.18)(38.75,26.13)(38.81,26.11)(38.88,26.06)(38.95,26.02)(39.01,25.97)(39.08,25.93)(39.15,25.88)(39.19,25.84)
\path
(39.19,25.84)(39.26,25.79)(39.3,25.74)(39.36,25.68)(39.41,25.63)(39.44,25.58)(39.48,25.52)(39.52,25.47)(39.56,25.4)(39.59,25.34)
\path
(39.59,25.34)(39.62,25.29)(39.66,25.22)(39.69,25.16)(39.7,25.11)(39.73,25.04)(39.75,25.0)(39.76,24.93)(39.76,24.88)(39.77,24.81)
\path
(39.77,24.81)(39.77,24.75)(39.77,24.7)(39.77,24.63)(39.76,24.58)(39.76,24.52)(39.75,24.47)(39.73,24.4)(39.72,24.36)(39.69,24.29)
\path
(39.69,24.29)(39.66,24.25)(39.63,24.2)(39.61,24.15)(39.58,24.09)(39.54,24.04)(39.48,24.0)(39.44,23.95)(39.4,23.91)(39.33,23.88)
\path
(39.33,23.88)(39.29,23.84)(39.22,23.81)(39.16,23.77)(39.08,23.75)(39.01,23.72)(38.94,23.68)(38.86,23.66)(38.76,23.65)(38.68,23.63)
\path
(38.68,23.63)(38.58,23.61)(38.48,23.59)(38.38,23.59)(38.29,23.58)(38.18,23.58)(38.05,23.58)(37.94,23.59)(37.81,23.59)(37.69,23.59)
\path
(37.69,23.59)(37.55,23.61)(37.41,23.63)(37.27,23.65)(37.13,23.68)(36.98,23.72)(36.83,23.75)(36.68,23.79)(36.51,23.84)(36.34,23.88)
\path(36.34,23.88)(36.16,23.93)(36.0,24.0)(36.0,24.0)
\drawvector{37.13}{23.68}{1.13}{-1}{0}
\drawcenteredtext{23.4}{27.0}{$a/b_1,c;2$}
\drawcenteredtext{35.0}{20.0}{$a/b_2,c;1$}
\drawcenteredtext{37.3}{28.0}{$a/b_1,c;3$}
\drawdotline{26.0}{40.0}{22.0}{32.0}
\drawdotline{26.0}{16.0}{22.0}{8.0}
\drawdotline{28.0}{16.0}{30.0}{8.0}
\end{picture}
}

The following two examples have both sequential interfaces $\emptyset$ and
hence we will omit the sequential information.

\subsubsection{A philosopher}

Consider the alphabet $A=\left\{  t,r,\varepsilon\right\}  $. \ A philosopher
is an automaton $\mathbf{Phil}$ with left interface $A$ and right interfaces
$A$, state space $\left\{  1,2,3,4\right\}  $, both sequential interfaces
$\emptyset\subseteq\left\{  1,2,3,4\right\}  $, and transition matrices
\begin{align*}
\mathsf{Phil}_{\varepsilon,\varepsilon}  &  =\left[
\begin{array}
[c]{cccc}%
\frac{1}{2} & 0 & 0 & 0\\
0 & \frac{1}{2} & 0 & 0\\
0 & 0 & \frac{1}{2} & 0\\
0 & 0 & 0 & \frac{1}{2}%
\end{array}
\right]  ,\;\\
\mathsf{Phil}_{t,\varepsilon}  &  =\left[
\begin{array}
[c]{cccc}%
0 & \frac{1}{2} & 0 & 0\\
0 & 0 & 0 & 0\\
0 & 0 & 0 & 0\\
0 & 0 & 0 & 0
\end{array}
\right]  ,\;\mathsf{Phil}_{\varepsilon,t}=\left[
\begin{array}
[c]{cccc}%
0 & 0 & 0 & 0\\
0 & 0 & \frac{1}{2} & 0\\
0 & 0 & 0 & 0\\
0 & 0 & 0 & 0
\end{array}
\right] \\
\mathsf{Phil}_{r,\varepsilon}  &  =\left[
\begin{array}
[c]{cccc}%
0 & 0 & 0 & 0\\
0 & 0 & 0 & 0\\
0 & 0 & 0 & \frac{1}{2}\\
0 & 0 & 0 & 0
\end{array}
\right]  ,\;\mathsf{Phil}_{\varepsilon,r}=\left[
\begin{array}
[c]{cccc}%
0 & 0 & 0 & 0\\
0 & 0 & 0 & 0\\
0 & 0 & 0 & 0\\
\frac{1}{2} & 0 & 0 & 0
\end{array}
\right]  .
\end{align*}

The other four transition matrices are zero matrices.

Notice that the total matrix of $\mathbf{Phil}$ is
\[
\left[
\begin{array}
[c]{cccc}%
\frac{1}{2} & \frac{1}{2} & 0 & 0\\
0 & \frac{1}{2} & \frac{1}{2} & 0\\
0 & 0 & \frac{1}{2} & \frac{1}{2}\\
\frac{1}{2} & 0 & 0 & \frac{1}{2}%
\end{array}
\right]  ,\text{ }%
\]
which is clearly stochastic, so $\mathbf{Phil}$ is a Markov automaton.

The intention behind these matrices is as follows: in all states the
philosopher does a transition labelled $\varepsilon,\varepsilon$ (\emph{idle
transition}) with probability $\frac{1}{2}$; in state $1$ he does a transition
to state $2$ with probability $\frac{1}{2}$ labelled $t,\varepsilon$
(\emph{take the left fork}); in state $2$ he does a transition to state $3$
with probability $\frac{1}{2}$ labelled $\varepsilon,t$ (\emph{take the right
fork}); in state $3$ he does a transition to state $4$ with probability
$\frac{1}{2}$ labelled $r,\varepsilon$ (\emph{release the left fork}); and in
state $4$ he does a transition to state $1$ with probability $\frac{1}{2}$
labelled $\varepsilon,r$ (\emph{release the right fork}). All this information
may be put in the following diagram.

%
%
%
%
%
%
%
%
%
%
%
%
%
%
%
%
%
%
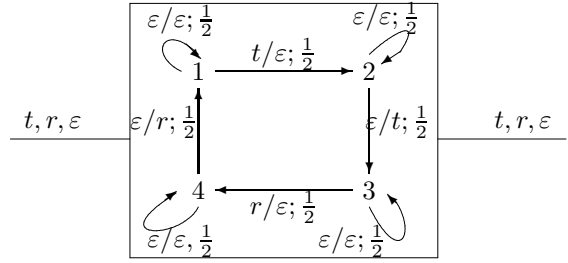
\centerline{\tt\setlength{\unitlength}{3.2pt}\begin{picture}(40,40)
\thinlines
\drawpath{4.0}{36.0}{40.0}{36.0}
\drawpath{40.0}{36.0}{40.0}{6.0}
\drawpath{40.0}{6.0}{4.0}{6.0}
\drawpath{4.0}{6.0}{4.0}{36.0}
\drawpath{40.0}{20.0}{55.0}{20.0}
\drawpath{-10.0}{20.0}{4.0}{20.0}
\path
(10.0,28.0)(10.0,28.0)(9.88,28.05)(9.76,28.12)(9.64,28.18)(9.53,28.23)(9.43,28.3)(9.34,28.37)(9.23,28.44)(9.14,28.51)
\path
(9.14,28.51)(9.06,28.57)(8.97,28.63)(8.89,28.7)(8.81,28.78)(8.72,28.85)(8.65,28.93)(8.59,29.0)(8.52,29.06)(8.46,29.13)
\path
(8.46,29.13)(8.39,29.2)(8.34,29.28)(8.27,29.35)(8.22,29.43)(8.18,29.5)(8.14,29.56)(8.1,29.63)(8.06,29.7)(8.02,29.78)
\path
(8.02,29.78)(7.98,29.86)(7.96,29.93)(7.92,30.0)(7.9,30.06)(7.88,30.12)(7.86,30.2)(7.84,30.27)(7.84,30.34)(7.82,30.39)
\path
(7.82,30.39)(7.82,30.46)(7.82,30.53)(7.8,30.59)(7.82,30.64)(7.82,30.7)(7.82,30.77)(7.82,30.82)(7.84,30.88)(7.86,30.94)
\path
(7.86,30.94)(7.88,31.0)(7.88,31.04)(7.92,31.1)(7.94,31.14)(7.96,31.2)(7.98,31.25)(8.02,31.29)(8.06,31.32)(8.1,31.37)
\path
(8.1,31.37)(8.13,31.4)(8.17,31.45)(8.22,31.47)(8.26,31.51)(8.31,31.54)(8.35,31.56)(8.4,31.59)(8.46,31.62)(8.52,31.63)
\path
(8.52,31.63)(8.57,31.65)(8.64,31.67)(8.69,31.68)(8.77,31.69)(8.82,31.7)(8.89,31.7)(8.97,31.7)(9.03,31.7)(9.11,31.7)
\path
(9.11,31.7)(9.18,31.69)(9.27,31.68)(9.35,31.67)(9.43,31.64)(9.52,31.62)(9.6,31.61)(9.68,31.57)(9.77,31.54)(9.86,31.52)
\path
(9.86,31.52)(9.96,31.47)(10.05,31.44)(10.14,31.38)(10.23,31.34)(10.34,31.29)(10.43,31.22)(10.55,31.17)(10.64,31.11)(10.75,31.04)
\path
(10.75,31.04)(10.85,30.96)(10.97,30.88)(11.07,30.8)(11.18,30.71)(11.3,30.62)(11.4,30.54)(11.52,30.44)(11.64,30.32)(11.76,30.22)
\path(11.76,30.22)(11.88,30.11)(11.98,30.0)(12.0,30.0)
\path
(32.0,30.0)(32.0,30.0)(32.04,30.05)(32.11,30.11)(32.18,30.18)(32.25,30.23)(32.31,30.3)(32.38,30.36)(32.47,30.43)(32.54,30.5)
\path
(32.54,30.5)(32.61,30.55)(32.68,30.62)(32.77,30.7)(32.86,30.77)(32.93,30.84)(33.02,30.89)(33.11,30.96)(33.18,31.04)(33.27,31.11)
\path
(33.27,31.11)(33.36,31.18)(33.45,31.25)(33.54,31.3)(33.63,31.37)(33.72,31.45)(33.81,31.52)(33.9,31.59)(34.0,31.64)(34.08,31.7)
\path
(34.08,31.7)(34.15,31.77)(34.25,31.84)(34.34,31.89)(34.43,31.95)(34.52,32.0)(34.61,32.06)(34.7,32.13)(34.79,32.18)(34.86,32.22)
\path
(34.86,32.22)(34.95,32.27)(35.04,32.33)(35.11,32.36)(35.2,32.4)(35.29,32.45)(35.36,32.5)(35.43,32.54)(35.5,32.56)(35.59,32.59)
\path
(35.59,32.59)(35.65,32.63)(35.72,32.65)(35.79,32.68)(35.86,32.7)(35.93,32.72)(36.0,32.75)(36.04,32.75)(36.11,32.77)(36.15,32.77)
\path
(36.15,32.77)(36.2,32.77)(36.25,32.77)(36.29,32.77)(36.34,32.77)(36.38,32.75)(36.43,32.75)(36.45,32.72)(36.47,32.7)(36.5,32.68)
\path
(36.5,32.68)(36.54,32.65)(36.54,32.61)(36.56,32.59)(36.58,32.54)(36.58,32.5)(36.58,32.45)(36.58,32.38)(36.58,32.33)(36.58,32.27)
\path
(36.58,32.27)(36.54,32.2)(36.54,32.13)(36.5,32.04)(36.5,31.95)(36.45,31.86)(36.43,31.77)(36.36,31.68)(36.33,31.57)(36.27,31.46)
\path
(36.27,31.46)(36.22,31.34)(36.15,31.21)(36.08,31.09)(36.0,30.96)(35.93,30.82)(35.84,30.68)(35.75,30.54)(35.65,30.37)(35.54,30.21)
\path
(12.0,12.0)(12.0,12.0)(11.93,11.88)(11.85,11.75)(11.78,11.64)(11.72,11.52)(11.63,11.43)(11.53,11.31)(11.47,11.19)(11.38,11.1)
\path
(11.38,11.1)(11.27,11.0)(11.18,10.9)(11.07,10.81)(10.97,10.72)(10.85,10.63)(10.75,10.56)(10.64,10.47)(10.53,10.38)(10.4,10.31)
\path
(10.4,10.31)(10.28,10.23)(10.18,10.14)(10.06,10.09)(9.93,10.02)(9.81,9.94)(9.68,9.88)(9.56,9.82)(9.43,9.77)(9.28,9.72)
\path
(9.28,9.72)(9.18,9.65)(9.03,9.6)(8.9,9.56)(8.78,9.52)(8.64,9.48)(8.52,9.44)(8.39,9.39)(8.27,9.38)(8.14,9.34)
\path
(8.14,9.34)(8.02,9.31)(7.9,9.28)(7.78,9.27)(7.65,9.25)(7.53,9.23)(7.42,9.22)(7.3,9.22)(7.19,9.19)(7.09,9.19)
\path
(7.09,9.19)(6.98,9.19)(6.88,9.19)(6.76,9.19)(6.67,9.22)(6.59,9.22)(6.5,9.25)(6.4,9.25)(6.3,9.27)(6.23,9.28)
\path
(6.23,9.28)(6.15,9.32)(6.09,9.35)(6.01,9.38)(5.94,9.4)(5.9,9.44)(5.84,9.5)(5.78,9.53)(5.75,9.59)(5.71,9.63)
\path
(5.71,9.63)(5.67,9.69)(5.65,9.75)(5.63,9.81)(5.63,9.88)(5.61,9.93)(5.61,10.0)(5.63,10.09)(5.63,10.14)(5.65,10.23)
\path
(5.65,10.23)(5.67,10.31)(5.71,10.39)(5.75,10.5)(5.8,10.59)(5.86,10.68)(5.92,10.77)(6.0,10.88)(6.07,10.98)(6.15,11.09)
\path
(6.15,11.09)(6.26,11.19)(6.38,11.32)(6.48,11.44)(6.61,11.56)(6.75,11.68)(6.88,11.81)(7.03,11.94)(7.19,12.09)(7.38,12.22)
\path
(7.38,12.22)(7.55,12.35)(7.75,12.52)(7.94,12.65)(8.18,12.81)(8.39,12.97)(8.63,13.13)(8.88,13.28)(9.14,13.47)(9.4,13.63)
\path
(34.83,13.0)(34.93,12.88)(35.0,12.75)(35.09,12.6)(35.18,12.5)(35.25,12.35)(35.33,12.23)(35.4,12.1)(35.47,11.98)(35.52,11.85)
\path
(35.52,11.85)(35.59,11.72)(35.65,11.6)(35.7,11.47)(35.75,11.35)(35.79,11.22)(35.84,11.1)(35.88,10.98)(35.93,10.85)(35.97,10.75)
\path
(35.97,10.75)(36.0,10.63)(36.02,10.52)(36.04,10.39)(36.08,10.28)(36.09,10.18)(36.11,10.06)(36.13,9.94)(36.15,9.84)(36.15,9.75)
\path
(36.15,9.75)(36.15,9.63)(36.15,9.53)(36.18,9.44)(36.15,9.34)(36.15,9.25)(36.15,9.15)(36.15,9.07)(36.13,9.0)(36.11,8.9)
\path
(36.11,8.9)(36.11,8.84)(36.09,8.75)(36.06,8.68)(36.04,8.6)(36.0,8.56)(36.0,8.5)(35.95,8.43)(35.93,8.38)(35.88,8.32)
\path
(35.88,8.32)(35.86,8.27)(35.81,8.23)(35.75,8.19)(35.72,8.15)(35.68,8.13)(35.61,8.1)(35.58,8.09)(35.52,8.06)(35.47,8.06)
\path
(35.47,8.06)(35.4,8.03)(35.34,8.03)(35.29,8.03)(35.22,8.06)(35.15,8.06)(35.08,8.09)(35.0,8.1)(34.93,8.13)(34.86,8.18)
\path
(34.86,8.18)(34.79,8.22)(34.72,8.25)(34.63,8.31)(34.54,8.35)(34.47,8.43)(34.38,8.5)(34.29,8.57)(34.2,8.64)(34.11,8.75)
\path
(34.11,8.75)(34.02,8.84)(33.93,8.93)(33.83,9.03)(33.75,9.14)(33.65,9.28)(33.54,9.39)(33.43,9.53)(33.33,9.68)(33.22,9.84)
\path
(33.22,9.84)(33.11,10.0)(33.0,10.14)(32.9,10.34)(32.79,10.52)(32.68,10.69)(32.58,10.89)(32.47,11.09)(32.34,11.31)(32.22,11.52)
\path(32.22,11.52)(32.11,11.75)(32.0,11.98)(32.0,12.0)
\drawvector{14.0}{28.0}{16.0}{1}{0}
\drawvector{32.0}{26.0}{10.0}{0}{-1}
\drawvector{30.0}{14.0}{16.0}{-1}{0}
\drawvector{12.0}{16.0}{10.0}{0}{1}
\drawvector{35.68}{30.37}{2.29}{-3}{-2}
\drawvector{8.07}{12.81}{1.39}{4}{3}
\drawvector{35.15}{12.46}{1.06}{-1}{1}
\drawvector{10.53}{31.25}{1.76}{1}{-1}
\drawcenteredtext{12.0}{28.0}{$1$}
\drawcenteredtext{32.0}{28.0}{$2$}
\drawcenteredtext{32.0}{14.0}{$3$}
\drawcenteredtext{12.0}{14.0}{$4$}
\drawcenteredtext{10.0}{34.0}{$\varepsilon/\varepsilon;\frac{1}{2}$}
\drawcenteredtext{34.0}{34.3}{$\varepsilon/\varepsilon;\frac{1}{2}$}
\drawcenteredtext{30.0}{8.0}{$\varepsilon/\varepsilon;\frac{1}{2}$}
\drawcenteredtext{10.0}{8.0}{$\varepsilon/\varepsilon,\frac{1}{2}$}
\drawcenteredtext{22.0}{30.0}{$t/\varepsilon;\frac{1}{2}$}
\drawcenteredtext{35.3}{22.0}{$\varepsilon/t;\frac{1}{2}$}
\drawcenteredtext{22.0}{12.0}{$r/\varepsilon;\frac{1}{2}$}
\drawcenteredtext{8.0}{22.0}{$\varepsilon/r;\frac{1}{2}$}
\drawcenteredtext{50.0}{22.0}{$t,r,\varepsilon$}
\drawcenteredtext{-5.0}{22.0}{$t,r,\varepsilon$}
\end{picture}}

\subsubsection{A fork}

Consider again the alphabet $A=\left\{  t,r,\varepsilon\right\}  $. \ A fork
is an automaton $\mathbf{Fork}$ with left interface $A$ and right interface
$A$, state space $\left\{  1,2,3\right\}  $, both sequential interfaces
$\emptyset\subseteq\left\{  1,2,3\right\}  $, and transition matrices
\begin{align*}
\mathsf{Fork}_{\varepsilon,\varepsilon}  &  =\left[
\begin{array}
[c]{ccc}%
\frac{1}{3} & 0 & 0\\
0 & \frac{1}{2} & 0\\
0 & 0 & \frac{1}{2}%
\end{array}
\right]  ,\;\\
\mathsf{Fork}_{t,\varepsilon}  &  =\left[
\begin{array}
[c]{ccc}%
0 & \frac{1}{3} & 0\\
0 & 0 & 0\\
0 & 0 & 0
\end{array}
\right]  ,\;\mathsf{Fork}_{\varepsilon,t}=\left[
\begin{array}
[c]{ccc}%
0 & 0 & \frac{1}{3}\\
0 & 0 & 0\\
0 & 0 & 0
\end{array}
\right] \\
\mathsf{Fork}_{r,\varepsilon}  &  =\left[
\begin{array}
[c]{ccc}%
0 & 0 & 0\\
\frac{1}{2} & 0 & 0\\
0 & 0 & 0
\end{array}
\right]  ,\;\mathsf{Fork}_{\varepsilon,r}=\left[
\begin{array}
[c]{ccc}%
0 & 0 & 0\\
0 & 0 & 0\\
\frac{1}{2} & 0 & 0
\end{array}
\right]  .
\end{align*}
The other four transition matrices are zero.

$\mathbf{Fork}$ is a Markov automaton since its total matrix is
\[
\left[
\begin{array}
[c]{ccc}%
\frac{1}{3} & \frac{1}{3} & \frac{1}{3}\\
\frac{1}{2} & \frac{1}{2} & 0\\
\frac{1}{2} & 0 & \frac{1}{2}%
\end{array}
\right]  .
\]

The intention behind these matrices is as follows: in all states the fork does
a transition labelled $\varepsilon,\varepsilon$ (\emph{idle transition}) with
positive probability (either $\frac{1}{3}$ or $\frac{1}{2}$); in state $1$ it
does a transition to state $2$ with probability $\frac{1}{3}$ labelled
$t,\varepsilon$ (\emph{taken to the left}); in state $1$ he does a transition
to state $3$ with probability $\frac{1}{3}$ labelled $\varepsilon,t$
(\emph{taken to the right}); in state $2$ he does a transition to state $1$
with probability $\frac{1}{2}$ labelled $r,\varepsilon$ (\emph{released to the
left}); in state $3$ he does a transition to state $1$ with probability
$\frac{1}{2}$ labelled $\varepsilon,r$ (\emph{released to the right}).

All this information may be put in the following diagram:

%
%
%
%
%
%
%
%
%
%
%
%
%
%
%
%
%
%
%
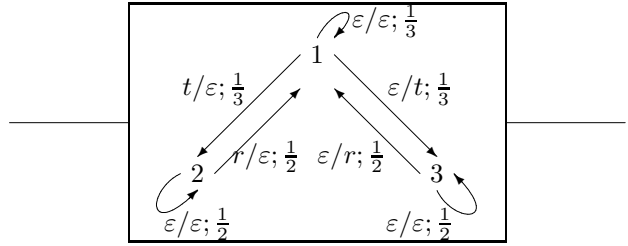
\centerline{\tt\setlength{\unitlength}{3.2pt}\begin{picture}(90,36)
\thinlines
\drawframebox{46.0}{18.0}{44.0}{28.0}{}
\drawpath{24.0}{18.0}{10.0}{18.0}
\drawpath{68}{18.0}{82.0}{18.0}
\drawcenteredtext{46.0}{26.0}{$1$}
\drawcenteredtext{32.0}{12.0}{$2$}
\drawcenteredtext{60.0}{12.0}{$3$}
\drawvector{48.0}{26.0}{12.0}{1}{-1}
\drawvector{58.0}{12.0}{10.0}{-1}{1}
\drawvector{44.0}{26.0}{12.0}{-1}{-1}
\drawvector{34.0}{12.0}{10.0}{1}{1}
\path
(30.0,12.0)(30.0,12.0)(29.88,11.93)(29.75,11.86)(29.65,11.79)(29.52,11.74)(29.43,11.66)(29.31,11.59)(29.2,11.52)(29.11,11.45)
\path
(29.11,11.45)(29.0,11.36)(28.91,11.29)(28.81,11.2)(28.72,11.13)(28.63,11.04)(28.56,10.95)(28.47,10.88)(28.38,10.79)(28.31,10.7)
\path
(28.31,10.7)(28.24,10.61)(28.15,10.52)(28.09,10.43)(28.02,10.34)(27.95,10.25)(27.88,10.16)(27.83,10.08)(27.77,10.0)(27.72,9.9)
\path
(27.72,9.9)(27.66,9.81)(27.61,9.72)(27.56,9.63)(27.52,9.54)(27.49,9.45)(27.45,9.36)(27.4,9.27)(27.38,9.2)(27.34,9.11)
\path
(27.34,9.11)(27.31,9.02)(27.29,8.93)(27.27,8.86)(27.25,8.77)(27.24,8.7)(27.22,8.61)(27.22,8.54)(27.2,8.47)(27.2,8.38)
\path
(27.2,8.38)(27.2,8.31)(27.2,8.25)(27.2,8.18)(27.22,8.11)(27.22,8.06)(27.25,8.0)(27.25,7.93)(27.27,7.88)(27.29,7.83)
\path
(27.29,7.83)(27.33,7.77)(27.36,7.72)(27.38,7.68)(27.41,7.63)(27.45,7.59)(27.5,7.56)(27.54,7.52)(27.59,7.5)(27.63,7.47)
\path
(27.63,7.47)(27.7,7.45)(27.75,7.43)(27.81,7.41)(27.88,7.4)(27.93,7.4)(28.0,7.4)(28.09,7.4)(28.15,7.4)(28.24,7.4)
\path
(28.24,7.4)(28.31,7.43)(28.4,7.43)(28.5,7.47)(28.59,7.49)(28.68,7.52)(28.77,7.56)(28.88,7.61)(28.99,7.65)(29.09,7.7)
\path
(29.09,7.7)(29.2,7.75)(29.33,7.83)(29.45,7.9)(29.56,7.97)(29.68,8.04)(29.81,8.13)(29.95,8.22)(30.09,8.33)(30.22,8.43)
\path
(30.22,8.43)(30.36,8.54)(30.52,8.65)(30.66,8.77)(30.81,8.9)(30.97,9.04)(31.13,9.18)(31.29,9.33)(31.47,9.49)(31.63,9.65)
\path(31.63,9.65)(31.81,9.81)(31.99,9.99)(32.0,10.0)
\drawvector{30.97}{9.04}{1.02}{1}{1}
\path
(60.0,10.0)(60.0,10.0)(60.05,9.88)(60.12,9.75)(60.19,9.65)(60.25,9.52)(60.31,9.43)(60.38,9.31)(60.47,9.2)(60.54,9.11)
\path
(60.54,9.11)(60.62,9.0)(60.69,8.91)(60.77,8.81)(60.86,8.72)(60.94,8.63)(61.02,8.56)(61.11,8.47)(61.19,8.38)(61.27,8.31)
\path
(61.27,8.31)(61.37,8.24)(61.45,8.15)(61.55,8.09)(61.63,8.02)(61.73,7.95)(61.81,7.88)(61.91,7.83)(62.0,7.77)(62.08,7.72)
\path
(62.08,7.72)(62.16,7.66)(62.26,7.61)(62.34,7.56)(62.44,7.52)(62.52,7.49)(62.62,7.45)(62.7,7.4)(62.79,7.38)(62.87,7.34)
\path
(62.87,7.34)(62.95,7.31)(63.05,7.29)(63.12,7.27)(63.2,7.25)(63.29,7.24)(63.37,7.22)(63.44,7.22)(63.51,7.2)(63.59,7.2)
\path
(63.59,7.2)(63.66,7.2)(63.73,7.2)(63.8,7.2)(63.87,7.22)(63.93,7.22)(64.0,7.25)(64.05,7.25)(64.11,7.27)(64.16,7.29)
\path
(64.16,7.29)(64.2,7.33)(64.26,7.36)(64.3,7.38)(64.34,7.41)(64.38,7.45)(64.43,7.5)(64.45,7.54)(64.48,7.59)(64.51,7.63)
\path
(64.51,7.63)(64.54,7.7)(64.55,7.75)(64.56,7.81)(64.58,7.88)(64.58,7.93)(64.58,8.0)(64.58,8.09)(64.58,8.15)(64.58,8.24)
\path
(64.58,8.24)(64.55,8.31)(64.55,8.4)(64.51,8.5)(64.5,8.59)(64.45,8.68)(64.43,8.77)(64.37,8.88)(64.33,8.99)(64.27,9.09)
\path
(64.27,9.09)(64.23,9.2)(64.16,9.33)(64.08,9.45)(64.01,9.56)(63.94,9.68)(63.84,9.81)(63.76,9.95)(63.66,10.09)(63.55,10.22)
\path
(63.55,10.22)(63.44,10.36)(63.33,10.52)(63.2,10.66)(63.08,10.81)(62.94,10.97)(62.8,11.13)(62.66,11.29)(62.5,11.47)(62.33,11.63)
\path(62.33,11.63)(62.16,11.81)(62.0,11.99)(62.0,12.0)
\drawvector{62.94}{10.97}{0.94}{-1}{1}
\path
(46.0,28.0)(46.0,28.0)(46.05,28.11)(46.12,28.22)(46.18,28.34)(46.23,28.45)(46.3,28.56)(46.37,28.66)(46.44,28.77)(46.51,28.88)
\path
(46.51,28.88)(46.58,28.97)(46.65,29.08)(46.7,29.16)(46.77,29.25)(46.84,29.34)(46.93,29.43)(47.0,29.52)(47.06,29.61)(47.13,29.68)
\path
(47.13,29.68)(47.2,29.77)(47.27,29.84)(47.36,29.9)(47.43,29.99)(47.5,30.04)(47.56,30.11)(47.63,30.18)(47.7,30.25)(47.77,30.29)
\path
(47.77,30.29)(47.86,30.36)(47.93,30.4)(48.0,30.47)(48.05,30.52)(48.12,30.56)(48.19,30.61)(48.26,30.65)(48.33,30.68)(48.4,30.72)
\path
(48.4,30.72)(48.47,30.75)(48.52,30.79)(48.58,30.81)(48.65,30.84)(48.72,30.86)(48.76,30.9)(48.83,30.91)(48.88,30.93)(48.94,30.95)
\path
(48.94,30.95)(49.0,30.97)(49.05,30.97)(49.09,30.97)(49.15,30.99)(49.19,30.99)(49.25,31.0)(49.29,30.99)(49.33,30.99)(49.37,30.97)
\path
(49.37,30.97)(49.41,30.97)(49.44,30.97)(49.48,30.95)(49.51,30.93)(49.54,30.91)(49.56,30.9)(49.58,30.88)(49.62,30.84)(49.63,30.81)
\path
(49.63,30.81)(49.66,30.79)(49.66,30.75)(49.68,30.72)(49.69,30.68)(49.69,30.65)(49.7,30.61)(49.7,30.56)(49.69,30.52)(49.69,30.47)
\path
(49.69,30.47)(49.69,30.4)(49.68,30.36)(49.66,30.29)(49.65,30.25)(49.62,30.18)(49.61,30.11)(49.58,30.04)(49.55,29.99)(49.51,29.91)
\path
(49.51,29.91)(49.48,29.84)(49.44,29.77)(49.38,29.68)(49.33,29.61)(49.29,29.52)(49.23,29.43)(49.16,29.34)(49.11,29.25)(49.04,29.16)
\path
(49.04,29.16)(48.97,29.08)(48.88,28.97)(48.8,28.88)(48.72,28.77)(48.62,28.66)(48.54,28.56)(48.44,28.45)(48.33,28.34)(48.23,28.22)
\path(48.23,28.22)(48.11,28.11)(48.0,28.0)(48.0,28.0)
\drawvector{48.88}{28.97}{0.88}{-1}{-1}
\drawcenteredtext{54.0}{30.0}{$\varepsilon/\varepsilon;\frac{1}{3}$}
\drawcenteredtext{32.0}{6.0}{$\varepsilon/\varepsilon;\frac{1}{2}$}
\drawcenteredtext{58.0}{6.0}{$\varepsilon/\varepsilon;\frac{1}{2}$}
\drawcenteredtext{34.0}{22.0}{$t/\varepsilon;\frac{1}{3}$}
\drawcenteredtext{58.0}{22.0}{$\varepsilon/t;\frac{1}{3}$}
\drawcenteredtext{40.0}{14.0}{$r/\varepsilon;\frac{1}{2}$}
\drawcenteredtext{50.0}{14.0}{$\varepsilon/r;\frac{1}{2}$}
\end{picture}
}

\section{The algebra of weighted automata:\ operations}

Now we define operations on weighted automata analogous (in a precise sense)to
those defined in \cite{KSW97b, KSW00}.

\subsection{Sequential operations}

\subsubsection{Sum}

\begin{dfn}
Given weighted automata $\mathbf{Q}_{Y;A,B}^{X}$ and $\mathbf{R}_{W;C,D}^{Z}$
the \emph{sum} $\mathbf{Q\boxplus R}$ is the weighted automaton which has set
of states the disjoint union $Q+R$, left interfaces $A+C$, right interface
$B+D,$ top interface $X+Z$, bottom interface $Y+W$, (all disjoint sums),
$\gamma_{0}=\gamma_{0,\mathbf{Q}}+\gamma_{0,\mathbf{R}}$, $\gamma_{1}%
=\gamma_{1,\mathbf{Q}}+\gamma_{1,\mathbf{R}}$. The transition matrices are
\begin{align*}
\lbrack(\mathsf{Q+R})_{a,b}]_{q,q^{\prime}} &  =[\mathsf{Q}_{a,b}%
]_{q,q^{\prime}},\\
\lbrack(\mathsf{Q+R})_{c,d}]_{r,r^{\prime}} &  =[\mathsf{R}_{c,d}%
]_{r,r^{\prime}},
\end{align*}
all other values being $0$.
\end{dfn}

\subsubsection{Sequential composition}

\begin{dfn}
Given weighted automata $\mathbf{Q}_{Y;A,B}^{X}$ and $\mathbf{R}_{Z;C,D}^{Y}$,
the \emph{sequential composite of weighted automata} $\mathbf{Q\circ
}\mathbf{R}$ has set of states the equivalence classes of $Q+R$ under the
equivalence relation generated by the relation $\gamma_{1,\mathbf{Q}}%
(y)\sim\gamma_{0,\mathbf{R}}(y)$, $(y\in Y)$. The left interface is the
disjoint sum $A+C$, right interface $B+D$, the top interface is $X$, the
bottom interface is $Z$. The interface functions are
\[
\gamma_{0}=X\overset{\gamma_{0}}{\rightarrow}Q\rightarrow Q+R\rightarrow
(Q+R)/\sim\text{, }\gamma_{1}=Z\overset{\gamma_{1}}{\rightarrow}R\rightarrow
Q+R\rightarrow(Q+R)/\sim.
\]
Denoting the equivalence class of a state $s$ by $[s]$ the transition matrices
are:
\begin{align*}
\lbrack(\mathsf{Q\circ R})_{a,b}]_{[q],[q^{\prime}]} &  =\sum_{s\in\lbrack
q],s^{\prime}\in\lbrack q^{\prime}]}[\mathsf{Q}_{a,b}]_{s,s^{\prime}},\\
\lbrack(\mathsf{Q\circ R})_{c,d}]_{[r],[r^{\prime}]} &  =\sum_{s\in\lbrack
r],s^{\prime}\in\lbrack r^{\prime}]}[\mathsf{R}_{c,d}]_{s,s^{\prime}},
\end{align*}
\newline all other values being $0$.
\end{dfn}

\subsubsection{Sequential constants}

\begin{dfn}
Given a relation $\rho\subset X+Y$ we define a weighted automaton
$\mathbf{Seq(\rho)}$ as follows: it has state space the equivalence classes of
$X+Y$ under the equivalence relation $\sim\ $generated by $\rho$. It has
parallel interfaces $\emptyset$, so there are no transition matrices. The
sequential interfaces are $\gamma_{0}:X\rightarrow(X+Y)/\sim$, $\gamma
_{1}:Y\rightarrow(X+Y)/\sim$, both functions taking an element to its
equivalence class.
\end{dfn}

\paragraph{Sequential connectors}

Some special cases have particular importance and are called \emph{sequential
connectors} or \emph{wires}: (i) the automaton corresponding to the identity
function $1_{A}$, considered as a relation on $A\times A$ is also called
$1_{A}$; (ii) the automaton corresponding to the codiagonal function
$\nabla:A+A\rightarrow A$ (considered as a relation)\thinspace is called
$\nabla$; the automaton corresponding to the opposite relation $\nabla$ is
called $\nabla^{o}$; (iii) the automaton corresponding to the function
$twist:A\times B\rightarrow B\times A$ is called $twist_{A,B}$; (iv) the
automaton corresponding to the function $\emptyset\subseteq A$ is called $i$;
the automaton corresponding to the opposite relation of the function
$\emptyset\subseteq A$ is called $i^{o}$. Notice that we have overworked the
symbol $\nabla$ and it will be used again in another sense; however the
context should make clear which use we have in mind.

The role of sequential wires is to\emph{ equate states}.

\paragraph{The distributive law}

The bijection $\delta:X\times Y+X\times Z\rightarrow X\times(Y+Z)$ and its
inverse $\delta^{-1}:X\times(Y+Z)\rightarrow X\times Y+X\times Z)$ considered
as relations yield weighted automata which we will refer to with the same names.

\subsection{Parallel operations}

\subsubsection{Parallel composition}

\begin{dfn}
Given weighted automata $\mathbf{Q}_{Y;A,B}^{X}$ and $\mathbf{R}_{W;C,D}^{Z}$
the \emph{parallel composite} $\mathbf{Q\times R}$ is the weighted automaton
which has set of states $Q\times R$, left interfaces $A\times C$, right
interface $B\times D$, top interface $X\times Z$, bottom interface $Y\times
W$, sequential interface functions $\gamma_{0,\mathbf{Q}}\times\gamma
_{0,\mathbf{R}}$, $\gamma_{1,\mathbf{Q}}\times\gamma_{1,\mathbf{R}}$ and
transition matrices$,$
\[
(\mathsf{Q\times R})_{(a,c),(b,d)}=\mathsf{Q}_{a,b}\otimes\mathsf{R}_{c,d}.
\]
If the automata are positive weighted then $\varepsilon_{A\times
C}=(\varepsilon_{A},\varepsilon_{C}),$ $\varepsilon_{B\times D}=(\varepsilon
_{B},\varepsilon_{D})$.
\end{dfn}

This just says that the weight of a transition from $(q,r)$ to $(q^{\prime
},r^{\prime})$ with left signal $(a,c)$ and right signal $(b,d)$ is the
product of the weights of the transition $q\rightarrow q^{\prime}$ with
signals $a$ and $b$, and the weight of the transition $r\rightarrow r^{\prime
}$ with signals $c$ and $d$. The following simple lemma was proved in
\cite{dFSW09a}.

\begin{lem}
If $\mathbf{Q}$ and $\mathbf{R}$ are positive weighted automata then so is
$\mathbf{Q\times R}$ and
\[
\mathbf{N(Q\times R)=N(Q)\times N(R)}.
\]
Hence if $\mathbf{Q}$ and $\mathbf{R}$ are Markov automata then so is
$\mathbf{Q\times R}$.
\end{lem}

\begin{lem}
If $\mathbf{Q}$ and $\mathbf{R}$ are Markov automata then $(\mathbf{Q}%
\times\mathbf{R)}^{k}=\mathbf{Q}^{k}\times\mathbf{R}^{k}.$
\end{lem}

\subsubsection{Parallel with communication}

\begin{dfn}
Given weighted automata $\mathbf{Q}_{Y;A,B}^{X}$ and $\mathbf{R}_{W;B,C}^{Z}$
the\emph{ communicating parallel composite of weighted automata}
$\mathbf{Q||}\mathbf{R}$ has set of states $Q\times R$, left interfaces $A$,
right interface $C$,in and out interfaces $X\times Z$, $Y\times W$, sequential
interfaces $\gamma_{0,\mathbf{Q}}\times\gamma_{0,\mathbf{R}}$, $\gamma
_{1,\mathbf{Q}}\times\gamma_{1,\mathbf{R}}$ and transition matrices
\[
(\mathsf{Q||R})_{a,c}=\sum_{b\in B}\mathsf{Q}_{a,b}\otimes\mathsf{R}_{b,c}.
\]
\end{dfn}

The following simple lemmas were proved in \cite{dFSW09a}.

\begin{lem}
If $\mathbf{Q}$ and $\mathbf{R}$ are positive weighted automata then so is
$\mathbf{Q||R}$ $\ $and $\mathbf{N(N(Q)||N(}\mathbf{R))=N(Q||R).}$
\end{lem}

\subsubsection{Parallel constants}

\begin{dfn}
Given a relation $\rho\subseteq A\times B$ we define a weighted automaton
$\mathbf{Par}(\rho)$ as follows: it has one state $\ast$ say. Top and bottom
interfaces have one element. The transition matrices $[\mathsf{Par(\rho
)}_{a,b}]$ are $1\times1$ matrices, that is, real numbers. Then
$\mathsf{Par(\rho)}_{a,b}=1$ if $\rho$ relates $a$ and $b$, and
$\mathsf{Par(\rho)}_{a,b}=0$ otherwise. If $(\varepsilon_{A},\varepsilon
_{B})\in\rho$ then $Par(\rho)$ is also positive weighted.
\end{dfn}

\paragraph{Parallel connectors}

Some special cases, all described in \cite{KSW97b}, have particular importance
and are called \emph{parallel connectors }or \emph{wires}: (i) the automaton
corresponding to the identity function $1_{A}$, considered as a relation on
$A\times A$ is also called $1_{A}$; (ii) the automaton corresponding to the
diagonal function $\Delta:A\rightarrow A\times A$ (considered as a
relation)\thinspace is called $\Delta_{A}$; the automaton corresponding to the
opposite relation of $\Delta$ is called $\Delta_{A}^{o}$; (iii) the automaton
corresponding to the function $twist:A\times B\rightarrow B\times A$ is again
called $twist_{A,B}$; (iv) the automaton corresponding to the projection
function $A\rightarrow\{\ast\}$ is called $p$ and its opposite $p^{o}$.

The role of parallel wires is to \emph{equate signals}.

\paragraph{Parallel codiagonal}

The automaton corresponding to the function $\nabla:A+A\rightarrow A$ is
called the \emph{parallel codiagonal}, and is denoted $\nabla$ where is no
confusion. The automaton corresponding to the opposite relation is written
$\nabla^{o}$.

\subsection{Some derived operations}

\subsubsection{Local sum}

Given weighted automata $\mathbf{Q}_{Y;A,B}^{X}$ and $\mathbf{R}_{W;A,B}^{Z}$
the local sum $\mathbf{Q+R}$ is defined to be $\nabla_{A}^{o}%
||(\mathbf{Q\boxplus R)||}\nabla_{A}$. It has top and bottom interface $X+Z$
and $Y+W$, and left and right interfaces $A$ and $B$.

\subsubsection{ Local sequential composition}

Given weighted automata $\mathbf{Q}_{Y;A,B}^{X}$ and $\mathbf{R}_{Z;A,B}^{Y}$
the local sequential composite $\mathbf{Q\bullet R}$ is defined to be
$\nabla_{A}^{o}||(\mathbf{Q\circ R)||}\nabla_{A}$. It has top and bottom
interface $X$ and $Z$, and left and right interfaces $A$ and $B$.

\subsubsection{Sequential feedback}

Given weighted automata $\mathbf{Q}_{Y+Z;A,B}^{X+Z}$ sequential feedback
$\mathsf{Sfb}_{Z}(\mathbf{Q}_{Y+Z;A,B}^{X+Z})$ is defined to be
\[
(1_{X}\boxplus i_{Z})\circ(1_{X}\boxplus\nabla^{o}_{Z})\circ(\mathbf{Q}%
\boxplus1_{Z})\circ(1_{Y}\boxplus\nabla_{Z})\circ(1_{Y}\boxplus i_{Z}^{o}).
\]
This formula is easier to understand looking ahead at the graphical
representation in the next section.

\subsubsection{Parallel feedback}

Given weighted automata $\mathbf{Q}_{Y;A\times C,B\times C}^{X}$ parallel
feedback $\mathsf{Pfb}_{C}(\mathbf{Q}_{Y;A+C,B+C}^{X})$ is defined to be
\[
(1_{A}\times p_{C}^{o})||(1_{A}\times\Delta_{C})||(\mathbf{Q}\times
1_{C})||(1_{B}\times\Delta^{o}_{C})||(1_{B}\times p_{C}).
\]
This formula is also easier to understand looking ahead to the graphical
representation in the next section.

\subsection{Graphical representation of expressions of weighted automata}

Not only do weighted automata have a graphical representation, as seen above,
but so also do expressions of automata, as described in \cite{KSW97b}. We
extend the representation given in that paper to the combination of sequential
and parallel operations and constants. We will see in section 4.1 that the 
graphical representation of single automata
is actually a special case of this new representation of expressions, modulo
the equations satisfied by wires (the Frobenius and separable algebra equations
first introduces in \cite{CW87}, see also \cite{RSW05}).

In general an expression will be represented by a diagram of the following sort:

%
%
%
%
%
%
%
%
%
%
%
%
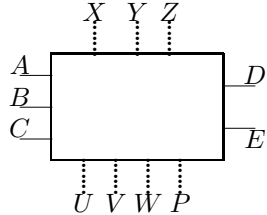
\centerline{\tt\setlength{\unitlength}{2pt}
\begin{picture}(54,48)
\thinlines
\drawframebox{26.0}{24.0}{32.0}{20.0}{}
\drawpath{10.0}{30.0}{4.0}{30.0}
\drawpath{10.0}{24.0}{4.0}{24.0}
\drawpath{10.0}{18.0}{4.0}{18.0}
\drawpath{42.0}{28.0}{48.0}{28.0}
\drawpath{42.0}{20.0}{48.0}{20.0}
\drawdotline{18.0}{40.0}{18.0}{34.0}
\drawdotline{26.0}{40.0}{26.0}{34.0}
\drawdotline{32.0}{40.0}{32.0}{34.0}
\drawdotline{16.0}{14.0}{16.0}{8.0}
\drawdotline{22.0}{14.0}{22.0}{8.0}
\drawdotline{28.0}{14.0}{28.0}{8.0}
\drawdotline{34.0}{14.0}{34.0}{8.0}
\drawcenteredtext{4.0}{32.0}{$A$}
\drawcenteredtext{4.0}{26.0}{$B$}
\drawcenteredtext{4.0}{20.0}{$C$}
\drawcenteredtext{48.0}{30.0}{$D$}
\drawcenteredtext{48.0}{18.0}{$E$}
\drawcenteredtext{18.0}{42.0}{$X$}
\drawcenteredtext{26.0}{42.0}{$Y$}
\drawcenteredtext{32.0}{42.0}{$Z$}
\drawcenteredtext{16.0}{6.0}{$U$}
\drawcenteredtext{22.0}{6.0}{$V$}
\drawcenteredtext{28.0}{6.0}{$W$}
\drawcenteredtext{34.0}{6.0}{$P$}
\end{picture}} The multiple lines on the left and right hand sides correspond
to parallel interfaces which are products of sets. For example, the component
has left interface $A\times B\times C$. Instead the multiple lines on the top
and bottom correspond to sequential interfaces which are disjoint sums of
sets, so the top interface is $X+Y+Z$.

\subsection{Operations and constants}

\subsubsection{Sum}

The sum $\mathbf{Q}^{X}_{Y;A,B}\boxplus\mathbf{R}^{Z}_{W;C,D}$ is represented as:

%
%
%
%
%
%
%
%
%
%
%
%
%
%
%
%
%
%
%
%
%
%
%
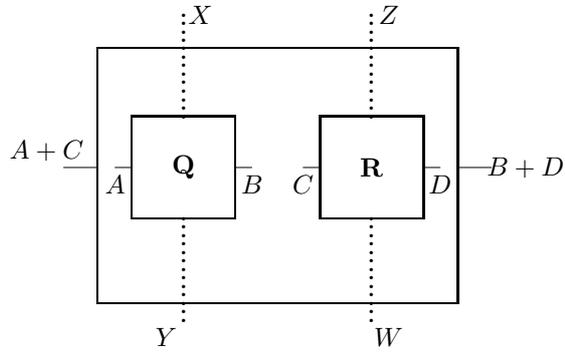
\centerline{\tt\setlength{\unitlength}{3.2pt}
\begin{picture}(70,52)
\thinlines
\drawframebox{22.0}{28.0}{12.0}{12.0}{$\mathbf{Q}$}
\drawpath{14.0}{28.0}{16.0}{28.0}
\drawpath{28.0}{28.0}{30.0}{28.0}
\drawframebox{44.0}{28.0}{12.0}{12.0}{$\mathbf{R}$}
\drawpath{36.0}{28.0}{38.0}{28.0}
\drawdotline{22.0}{34.0}{22.0}{46.0}
\drawdotline{44.0}{34.0}{44.0}{46.0}
\drawdotline{22.0}{22.0}{22.0}{10.0}
\drawdotline{44.0}{22.0}{44.0}{10.0}
\drawframebox{33.0}{27.0}{42.0}{30.0}{}
\drawcenteredtext{14.0}{26.0}{$A$}
\drawcenteredtext{30.0}{26.0}{$B$}
\drawpath{50.0}{28.0}{52.0}{28.0}
\drawcenteredtext{36.0}{26.0}{$C$}
\drawcenteredtext{52.0}{26.0}{$D$}
\drawcenteredtext{24.0}{46.0}{$X$}
\drawcenteredtext{20.0}{8.0}{$Y$}
\drawcenteredtext{46.0}{46.0}{$Z$}
\drawcenteredtext{46.0}{8.0}{$W$}
\drawpath{12.0}{28.0}{8.0}{28.0}
\drawpath{54.0}{28.0}{58.0}{28.0}
\drawcenteredtext{6.0}{30.0}{$A+C$}
\drawcenteredtext{62.0}{28.0}{$B+D$}
\end{picture}
}

\subsubsection{Sequential composition}

The sequential composite $\mathbf{Q}^{X}_{Y;A,B}\circ\mathbf{R}^{Y}_{Z;C,D}$
is represented as:

%
%
%
%
%
%
%
%
%
%
%
%
%
%
%
%
%
%
%
%
%
%
%
%
%
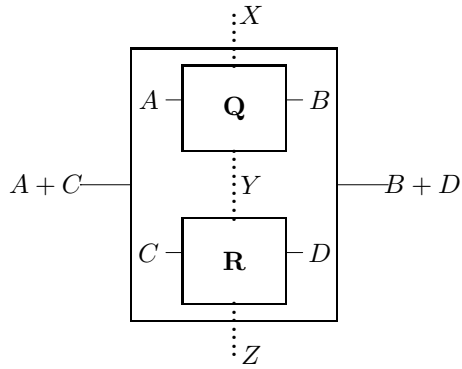
\centerline{\tt\setlength{\unitlength}{3.2pt}
\begin{picture}(58,58)
\thinlines
\drawframebox{28.0}{41.0}{12.0}{10.0}{$\mathbf{Q}$}
\drawframebox{28.0}{23.0}{12.0}{10.0}{$\mathbf{R}$}
\drawpath{20.0}{42.0}{22.0}{42.0}
\drawpath{34.0}{42.0}{36.0}{42.0}
\drawpath{20.0}{24.0}{22.0}{24.0}
\drawpath{34.0}{24.0}{36.0}{24.0}
\drawdotline{28.0}{52.0}{28.0}{46.0}
\drawdotline{28.0}{36.0}{28.0}{28.0}
\drawdotline{28.0}{18.0}{28.0}{12.0}
\drawframebox{28.0}{32.0}{24.0}{32.0}{}
\drawcenteredtext{30.0}{52.0}{$X$}
\drawcenteredtext{30.0}{32.0}{$Y$}
\drawcenteredtext{30.0}{12.0}{$Z$}
\drawpath{16.0}{32.0}{10.0}{32.0}
\drawpath{10.0}{32.0}{10.0}{32.0}
\drawpath{40.0}{32.0}{46.0}{32.0}
\drawcenteredtext{18.0}{42.0}{$A$}
\drawcenteredtext{38.0}{42.0}{$B$}
\drawcenteredtext{18.0}{24.0}{$C$}
\drawcenteredtext{38.0}{24.0}{$D$}
\drawcenteredtext{6.0}{32.0}{$A+C$}
\drawcenteredtext{50.0}{32.0}{$B+D$}
\end{picture}
}

\subsubsection{Sequential connectors}

The various sequential connectors are represented as:

%
%
%
%
%
%
%
%
%
%
%
%
%
%
%
%
%
%
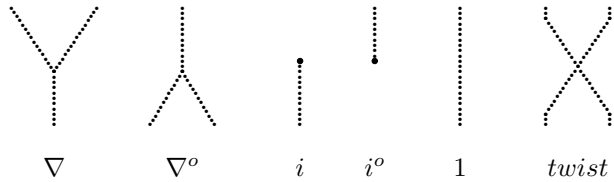
\centerline{\tt\setlength{\unitlength}{2pt}
\begin{picture}(120,40)
\thinlines
\drawdotline{4.0}{36.0}{12.0}{24.0}
\drawdotline{12.0}{24.0}{20.0}{36.0}
\drawdotline{12.0}{24.0}{12.0}{14.0}
\drawcenteredtext{12.0}{6.0}{$\nabla$}
\drawdotline{36.0}{36.0}{36.0}{24.0}
\drawdotline{36.0}{24.0}{30.0}{14.0}
\drawdotline{36.0}{24.0}{42.0}{14.0}
\drawcenteredtext{36.0}{6.0}{$\nabla^o$}
\drawdotline{72.0}{36.0}{72.0}{26.0}
\drawthickdot{72.0}{26.0}
\drawcenteredtext{58.0}{6.0}{$i$}
\drawdotline{58.0}{26.0}{58.0}{14.0}
\drawthickdot{58.0}{26.0}
\drawcenteredtext{72.0}{6.0}{$i^o$}
\drawdotline{88.0}{36.0}{88.0}{14.0}
\drawcenteredtext{88.0}{6.0}{$1$}
\drawdotline{104.0}{36.0}{104.0}{34.0}
\drawdotline{104.0}{34.0}{116.0}{16.0}
\drawdotline{116.0}{16.0}{116.0}{14.0}
\drawdotline{116.0}{36.0}{116.0}{34.0}
\drawdotline{116.0}{34.0}{104.0}{16.0}
\drawdotline{104.0}{16.0}{104.0}{14.0}
\drawcenteredtext{110.0}{6.0}{$twist$}
\end{picture}
}

\subsubsection{Distributive law}

The distributive law $\delta^{-1}:X\times(Y+Z)\to X\times Y+X\times Z$ and its
opposite are represented as:

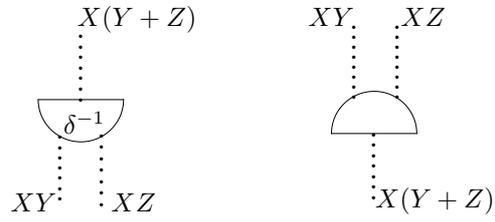
\centerline{\tt\setlength{\unitlength}{4pt}
\begin{picture}(52,28)
\thinlines
\drawdotline{36.12}{22.05}{36.12}{15.4}
\drawarc{38.0}{12.0}{8.0}{3.14}{6.28}
\drawpath{34.0}{12.0}{42.0}{12.0}
\drawdotline{37.97}{11.9}{37.97}{5.95}
\drawdotline{40.16}{22.05}{40.16}{15.4}
\drawcenteredtext{43.7}{5.7}{$X(Y+Z)$}
\drawcenteredtext{34.0}{23.0}{$XY$}
\drawcenteredtext{42.6}{23.0}{$XZ$}
\drawcenteredtext{38.0}{14.0}{}
\drawdotline{10.57}{21.2}{10.57}{15.2}
\drawarc{10.57}{15.2}{8.0}{6.28}{3.14}
\drawpath{6.57}{15.2}{14.57}{15.2}
\drawdotline{8.62}{11.62}{8.62}{5.86}
\drawdotline{12.49}{11.8}{12.49}{5.32}
\drawcenteredtext{10.91}{13.02}{$\delta^{-1}$}
\drawcenteredtext{15.83}{22.83}{$X(Y+Z)$}
\drawcenteredtext{6.16}{5.5}{$XY$}
\drawcenteredtext{15.48}{5.5}{$XZ$}
\end{picture}
}

\subsubsection{Product}

The product $\mathbf{Q}^{X}_{Y;A,B}\times\mathbf{R}^{Z}_{W;C,D}$ is
represented as:

%
%
%
%
%
%
%
%
%
%
%
%
%
%
%
%
%
%
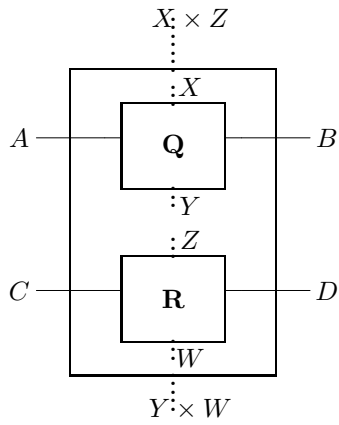
\centerline{\tt\setlength{\unitlength}{3.2pt}
\begin{picture}(46,64)
\thinlines
\drawframebox{22.0}{43.0}{12.0}{10.0}{$\mathbf{Q}$}
\drawframebox{22.0}{25.0}{12.0}{10.0}{$\mathbf{R}$}
\drawpath{14.0}{44.0}{16.0}{44.0}
\drawpath{28.0}{44.0}{30.0}{44.0}
\drawpath{14.0}{26.0}{16.0}{26.0}
\drawpath{28.0}{26.0}{30.0}{26.0}
\drawdotline{22.0}{50.0}{22.0}{48.0}
\drawdotline{22.0}{38.0}{22.0}{36.0}
\drawdotline{22.0}{32.0}{22.0}{30.0}
\drawframebox{22.0}{34.0}{24.0}{36.0}{}
\drawcenteredtext{4.0}{26.0}{$C$}
\drawcenteredtext{24.0}{36.0}{$Y$}
\drawcenteredtext{40.0}{26.0}{$D$}
\drawpath{6.0}{44.0}{14.0}{44.0}
\drawpath{6.0}{26.0}{14.0}{26.0}
\drawpath{30.0}{44.0}{38.0}{44.0}
\drawcenteredtext{4.0}{44.0}{$A$}
\drawcenteredtext{40.0}{44.0}{$B$}
\drawpath{30.0}{26.0}{38.0}{26.0}
\drawdotline{22.0}{20.0}{22.0}{18.0}
\drawdotline{22.0}{58.0}{22.0}{52.0}
\drawdotline{22.0}{16.0}{22.0}{12.0}
\drawcenteredtext{24.0}{12.0}{$Y\times W$}
\drawcenteredtext{24.0}{50.0}{$X$}
\drawcenteredtext{24.0}{32.0}{$Z$}
\drawcenteredtext{24.0}{18.0}{$W$}
\drawcenteredtext{24.0}{58.0}{$X\times Z$}
\end{picture}
}

\subsubsection{Parallel with communication}

The parallel with communication $\mathbf{Q}^{X}_{Y;A,B}||\mathbf{R}%
^{Z}_{W;B,C}$ is represented as:

%
%
%
%
%
%
%
%
%
%
%
%
%
%
%
%
%
%
%
%
%
%
%
%
%
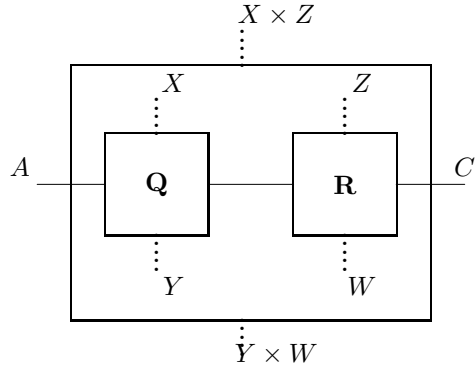
\centerline{\tt\setlength{\unitlength}{3.2pt}
\begin{picture}(62,58)
\thinlines
\drawframebox{20.0}{32.0}{12.0}{12.0}{$\mathbf{Q}$}
\drawpath{14.0}{32.0}{6.0}{32.0}
\drawpath{26.0}{32.0}{36.0}{32.0}
\drawframebox{42.0}{32.0}{12.0}{12.0}{$\mathbf{R}$}
\drawpath{48.0}{32.0}{56.0}{32.0}
\drawdotline{20.0}{42.0}{20.0}{38.0}
\drawdotline{20.0}{26.0}{20.0}{22.0}
\drawdotline{42.0}{42.0}{42.0}{38.0}
\drawdotline{42.0}{26.0}{42.0}{22.0}
\drawframebox{31.0}{31.0}{42.0}{30.0}{}
\drawdotline{30.0}{16.0}{30.0}{12.0}
\drawcenteredtext{4.0}{34.0}{$A$}
\drawcenteredtext{56.0}{34.0}{$C$}
\drawcenteredtext{22.0}{44.0}{$X$}
\drawcenteredtext{44.0}{44.0}{$Z$}
\drawcenteredtext{22.0}{20.0}{$Y$}
\drawcenteredtext{44.0}{20.0}{$W$}
\drawdotline{30.0}{46.0}{30.0}{50.0}
\drawcenteredtext{34.0}{52.0}{$X\times Z$}
\drawcenteredtext{34.0}{12.0}{$Y\times W$}
\end{picture}
}

\subsubsection{Parallel connectors}

The various parallel connectors are represented as:

%
%
%
%
%
%
%
%
%
%
%
%
%
%
%
%
%
%
%
%
%
%
%
%
%
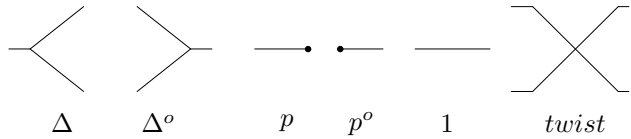
\centerline{\tt\setlength{\unitlength}{2.0pt}
\begin{picture}(124,32)
\thinlines
\drawpath{4.0}{20.0}{8.0}{20.0}
\drawpath{8.0}{20.0}{18.0}{28.0}
\drawpath{8.0}{20.0}{18.0}{12.0}
\drawpath{28.0}{28.0}{38.0}{20.0}
\drawpath{38.0}{20.0}{28.0}{12.0}
\drawpath{38.0}{20.0}{42.0}{20.0}
\drawpath{50.0}{20.0}{60.0}{20.0}
\drawthickdot{60.0}{20.0}
\drawpath{66.0}{20.0}{74.0}{20.0}
\drawthickdot{66.0}{20.0}
\drawpath{80.0}{20.0}{94.0}{20.0}
\drawpath{98.0}{28.0}{102.0}{28.0}
\drawpath{102.0}{28.0}{118.0}{12.0}
\drawpath{118.0}{12.0}{120.0}{12.0}
\drawpath{98.0}{12.0}{102.0}{12.0}
\drawpath{102.0}{12.0}{118.0}{28.0}
\drawpath{118.0}{28.0}{120.0}{28.0}
\drawcenteredtext{14.0}{6.0}{$\Delta$}
\drawcenteredtext{32.0}{6.0}{$\Delta^o$}
\drawcenteredtext{56.0}{6.0}{$p$}
\drawcenteredtext{70.0}{6.0}{$p^o$}
\drawcenteredtext{86.0}{6.0}{$1$}
\drawcenteredtext{110.0}{6.0}{$twist$}
\end{picture}
}

\subsubsection{Parallel codiagonal}

The parallel codiagonal $\nabla:A+A\to A$ and its opposite $\nabla^{o}$ are
represented as:

%
%
%
%
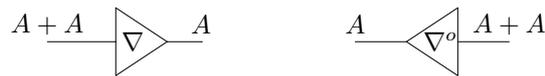
\centerline{\tt\setlength{\unitlength}{3.2pt}
\begin{picture}(68,16)
\thinlines
\drawpath{6.0}{8.0}{14.0}{8.0}
\drawpath{14.0}{8.0}{14.0}{12.0}
\drawpath{14.0}{12.0}{20.0}{8.0}
\drawpath{20.0}{8.0}{14.0}{4.0}
\drawpath{14.0}{4.0}{14.0}{8.0}
\drawpath{20.0}{8.0}{24.0}{8.0}
\drawpath{24.0}{8.0}{24.0}{8.0}
\drawcenteredtext{6.0}{10.0}{$A+A$}
\drawcenteredtext{24.0}{10.0}{$A$}
\drawcenteredtext{16.0}{8.0}{$\nabla$}
\drawpath{42.0}{8.0}{48.0}{8.0}
\drawpath{48.0}{8.0}{54.0}{12.0}
\drawpath{54.0}{12.0}{54.0}{4.0}
\drawpath{54.0}{4.0}{48.0}{8.0}
\drawpath{54.0}{8.0}{60.0}{8.0}
\drawcenteredtext{42.0}{10.0}{$A$}
\drawcenteredtext{60.0}{10.0}{$A+A$}
\drawcenteredtext{52.0}{8.0}{$\nabla^o$}
\end{picture}
}

\subsection{Some derived operations and constants}

\subsubsection{Repeated sequential and repeated parallel operations}

When representing repeated sequential operations, frequently we omit all but
the last bounding rectangle. With care this does not lead to ambiguity -
different interpretations lead to at worst isomorphic automata. We do the same
with repeated parallel operations. It is not possible to removed bounding
rectangles for mixed repeated and parallel operations without serious
ambiguity.

\subsubsection{Local Sum}

The local sum $\mathbf{Q}_{Y;A,B}^{X}+\mathbf{R}_{W;A,B}^{Z}$ is represented
by the first diagram below which includes the parallel codiagonals, but also,
since it is such a common derived operation, more briefly by the second
diagram below:

%
%
%
%
%
%
%
%
%
%
%
%
%
%
%
%
%
%
%
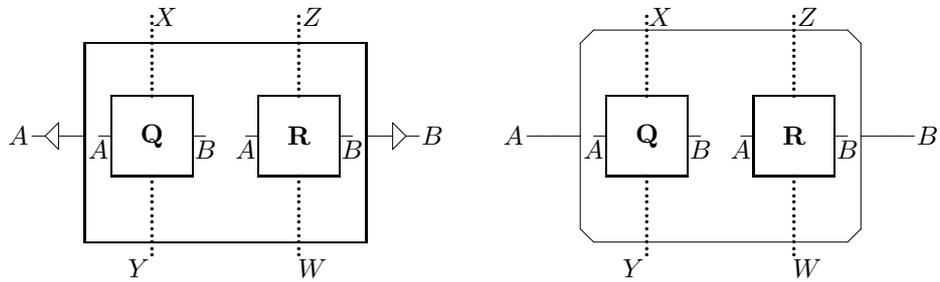
\centerline{\tt\setlength{\unitlength}{2.5pt}
\begin{picture}(72,52)
\thinlines
\drawframebox{24.0}{28.0}{12.0}{12.0}{$\mathbf{Q}$}
\drawpath{16.0}{28.0}{18.0}{28.0}
\drawpath{30.0}{28.0}{32.0}{28.0}
\drawframebox{46.0}{28.0}{12.0}{12.0}{$\mathbf{R}$}
\drawpath{38.0}{28.0}{40.0}{28.0}
\drawdotline{24.0}{34.0}{24.0}{46.0}
\drawdotline{46.0}{34.0}{46.0}{46.0}
\drawdotline{24.0}{22.0}{24.0}{10.0}
\drawdotline{46.0}{22.0}{46.0}{10.0}
\drawframebox{35.0}{27.0}{42.0}{30.0}{}
\drawcenteredtext{16.0}{26.0}{$A$}
\drawcenteredtext{32.0}{26.0}{$B$}
\drawpath{52.0}{28.0}{54.0}{28.0}
\drawcenteredtext{38.0}{26.0}{$A$}
\drawcenteredtext{54.0}{26.0}{$B$}
\drawcenteredtext{26.0}{46.0}{$X$}
\drawcenteredtext{22.0}{8.0}{$Y$}
\drawcenteredtext{48.0}{46.0}{$Z$}
\drawcenteredtext{48.0}{8.0}{$W$}
\drawpath{14.0}{28.0}{10.0}{28.0}
\drawpath{56.0}{28.0}{60.0}{28.0}
\drawpath{60.0}{30.0}{60.0}{26.0}
\drawpath{60.0}{26.0}{62.0}{28.0}
\drawpath{62.0}{28.0}{60.0}{30.0}
\drawpath{10.0}{30.0}{10.0}{26.0}
\drawpath{10.0}{26.0}{8.0}{28.0}
\drawpath{8.0}{28.0}{10.0}{30.0}
\drawpath{62.0}{28.0}{64.0}{28.0}
\drawpath{8.0}{28.0}{6.0}{28.0}
\drawcenteredtext{66.0}{28.0}{$B$}
\drawcenteredtext{4.0}{28.0}{$A$}
\end{picture}
\begin{picture}(72,52)
\thinlines
\drawframebox{24.0}{28.0}{12.0}{12.0}{$\mathbf{Q}$}
\drawpath{16.0}{28.0}{18.0}{28.0}
\drawpath{30.0}{28.0}{32.0}{28.0}
\drawframebox{46.0}{28.0}{12.0}{12.0}{$\mathbf{R}$}
\drawpath{38.0}{28.0}{40.0}{28.0}
\drawdotline{24.0}{34.0}{24.0}{46.0}
\drawdotline{46.0}{34.0}{46.0}{46.0}
\drawdotline{24.0}{22.0}{24.0}{10.0}
\drawdotline{46.0}{22.0}{46.0}{10.0}
\drawpath{14.0}{28.0}{14.0}{42.0}
\drawcenteredtext{16.0}{26.0}{$A$}
\drawcenteredtext{32.0}{26.0}{$B$}
\drawpath{52.0}{28.0}{54.0}{28.0}
\drawcenteredtext{38.0}{26.0}{$A$}
\drawcenteredtext{54.0}{26.0}{$B$}
\drawcenteredtext{26.0}{46.0}{$X$}
\drawcenteredtext{22.0}{8.0}{$Y$}
\drawcenteredtext{48.0}{46.0}{$Z$}
\drawcenteredtext{48.0}{8.0}{$W$}
\drawpath{14.0}{28.0}{10.0}{28.0}
\drawpath{56.0}{28.0}{60.0}{28.0}
\drawpath{8.0}{28.0}{10.0}{28.0}
\drawpath{60.0}{28.0}{62.0}{28.0}
\drawpath{14.0}{42.0}{16.0}{44.0}
\drawpath{16.0}{44.0}{54.0}{44.0}
\drawpath{54.0}{44.0}{56.0}{42.0}
\drawpath{56.0}{42.0}{56.0}{14.0}
\drawpath{62.0}{28.0}{64.0}{28.0}
\drawpath{8.0}{28.0}{6.0}{28.0}
\drawcenteredtext{66.0}{28.0}{$B$}
\drawcenteredtext{4.0}{28.0}{$A$}
\drawpath{56.0}{14.0}{54.0}{12.0}
\drawpath{54.0}{12.0}{16.0}{12.0}
\drawpath{16.0}{12.0}{14.0}{14.0}
\drawpath{14.0}{14.0}{14.0}{28.0}
\end{picture}
}

\subsubsection{Local sequential}

The local sequential $\mathbf{Q}_{Y;A,B}^{X}\bullet\mathbf{R}_{Z;A,B}^{Y}$ is
represented by the first diagram below which includes the parallel
codiagonals, but also, since it is such a common derived operation, more
briefly by the second diagram below:

%
%
%
%
%
%
%
%
%
%
%
%
%
%
%
%
%
%
%
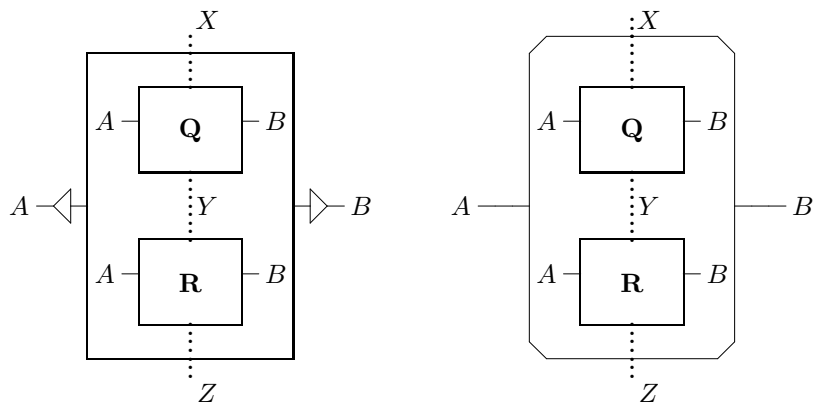
\centerline{\tt\setlength{\unitlength}{3.2pt}
\begin{picture}(50,60)
\thinlines
\drawframebox{24.0}{41.0}{12.0}{10.0}{$\mathbf{Q}$}
\drawframebox{24.0}{23.0}{12.0}{10.0}{$\mathbf{R}$}
\drawpath{16.0}{42.0}{18.0}{42.0}
\drawpath{30.0}{42.0}{32.0}{42.0}
\drawpath{16.0}{24.0}{18.0}{24.0}
\drawpath{30.0}{24.0}{32.0}{24.0}
\drawdotline{24.0}{52.0}{24.0}{46.0}
\drawdotline{24.0}{36.0}{24.0}{28.0}
\drawdotline{24.0}{18.0}{24.0}{12.0}
\drawframebox{24.0}{32.0}{24.0}{36.0}{}
\drawpath{12.0}{32.0}{10.0}{32.0}
\drawcenteredtext{26.0}{32.0}{$Y$}
\drawpath{10.0}{32.0}{10.0}{34.0}
\drawpath{10.0}{34.0}{8.0}{32.0}
\drawpath{8.0}{32.0}{10.0}{30.0}
\drawpath{10.0}{30.0}{10.0}{32.0}
\drawpath{8.0}{32.0}{6.0}{32.0}
\drawpath{36.0}{32.0}{38.0}{32.0}
\drawpath{38.0}{32.0}{38.0}{34.0}
\drawcenteredtext{26.0}{54.0}{$X$}
\drawcenteredtext{26.0}{10.0}{$Z$}
\drawpath{38.0}{34.0}{40.0}{32.0}
\drawpath{40.0}{32.0}{38.0}{30.0}
\drawpath{38.0}{30.0}{38.0}{32.0}
\drawpath{40.0}{32.0}{42.0}{32.0}
\drawcenteredtext{14.0}{42.0}{$A$}
\drawcenteredtext{14.0}{24.0}{$A$}
\drawcenteredtext{4.0}{32.0}{$A$}
\drawcenteredtext{34.0}{42.0}{$B$}
\drawcenteredtext{34.0}{24.0}{$B$}
\drawcenteredtext{44.0}{32.0}{$B$}
\end{picture}
\begin{picture}(50,60)
\thinlines
\drawframebox{24.0}{41.0}{12.0}{10.0}{$\mathbf{Q}$}
\drawframebox{24.0}{23.0}{12.0}{10.0}{$\mathbf{R}$}
\drawpath{16.0}{42.0}{18.0}{42.0}
\drawpath{30.0}{42.0}{32.0}{42.0}
\drawpath{16.0}{24.0}{18.0}{24.0}
\drawpath{30.0}{24.0}{32.0}{24.0}
\drawdotline{24.0}{52.0}{24.0}{46.0}
\drawdotline{24.0}{36.0}{24.0}{28.0}
\drawdotline{24.0}{18.0}{24.0}{12.0}
\drawpath{14.0}{52.0}{34.0}{52.0}
\drawpath{12.0}{32.0}{10.0}{32.0}
\drawcenteredtext{26.0}{32.0}{$Y$}
\drawpath{8.0}{32.0}{10.0}{32.0}
\drawpath{38.0}{32.0}{40.0}{32.0}
\drawdotline{24.0}{54.0}{24.0}{52.0}
\drawpath{8.0}{32.0}{6.0}{32.0}
\drawpath{36.0}{32.0}{38.0}{32.0}
\drawcenteredtext{26.0}{54.0}{$X$}
\drawcenteredtext{26.0}{10.0}{$Z$}
\drawpath{40.0}{32.0}{42.0}{32.0}
\drawcenteredtext{14.0}{42.0}{$A$}
\drawcenteredtext{14.0}{24.0}{$A$}
\drawcenteredtext{4.0}{32.0}{$A$}
\drawcenteredtext{34.0}{42.0}{$B$}
\drawcenteredtext{34.0}{24.0}{$B$}
\drawcenteredtext{44.0}{32.0}{$B$}
\drawpath{34.0}{52.0}{36.0}{50.0}
\drawpath{36.0}{50.0}{36.0}{16.0}
\drawpath{36.0}{16.0}{34.0}{14.0}
\drawpath{34.0}{14.0}{14.0}{14.0}
\drawpath{14.0}{14.0}{12.0}{16.0}
\drawpath{12.0}{16.0}{12.0}{50.0}
\drawpath{12.0}{50.0}{14.0}{52.0}
\end{picture}
}

\subsubsection{Sequential feedback}

The sequential feedback $\mathsf{Sfb}_{Z}(\mathbf{Q}_{Y+Z;A,B}^{X+Z})$ is
represented by the diagram:

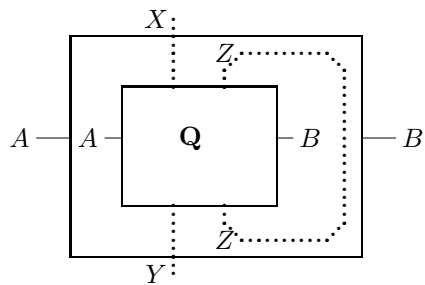
\centerline{\tt\setlength{\unitlength}{3.2pt}
\begin{picture}(56,42)
\thinlines
\drawframebox{25.0}{21.0}{18.0}{14.0}{}
\drawpath{16.0}{22.0}{14.0}{22.0}
\drawpath{34.0}{22.0}{36.0}{22.0}
\drawdotline{22.0}{32.0}{22.0}{30.0}
\drawdotline{22.0}{30.0}{22.0}{28.0}
\drawdotline{28.0}{30.0}{28.0}{28.0}
\drawdotline{28.0}{28.0}{28.0}{30.0}
\drawdotline{28.0}{30.0}{30.0}{32.0}
\drawdotline{30.0}{32.0}{40.0}{32.0}
\drawdotline{40.0}{32.0}{42.0}{30.0}
\drawdotline{42.0}{30.0}{42.0}{12.0}
\drawdotline{42.0}{12.0}{40.0}{10.0}
\drawdotline{40.0}{10.0}{32.0}{10.0}
\drawdotline{32.0}{10.0}{30.0}{10.0}
\drawdotline{30.0}{10.0}{28.0}{12.0}
\drawdotline{28.0}{12.0}{28.0}{14.0}
\drawdotline{22.0}{14.0}{22.0}{8.0}
\drawdotline{22.0}{32.0}{22.0}{34.0}
\drawcenteredtext{24.0}{22.0}{$\mathbf{Q}$}
\drawcenteredtext{12.0}{22.0}{$A$}
\drawcenteredtext{38.0}{22.0}{$B$}
\drawcenteredtext{20.0}{36.0}{$X$}
\drawcenteredtext{20.0}{6.0}{$Y$}
\drawcenteredtext{28.0}{32.0}{$Z$}
\drawcenteredtext{28.0}{10.0}{$Z$}
\drawframebox{27.0}{21.0}{34.0}{26.0}{}
\drawpath{10.0}{22.0}{6.0}{22.0}
\drawpath{44.0}{22.0}{48.0}{22.0}
\drawdotline{22.0}{34.0}{22.0}{36.0}
\drawdotline{22.0}{8.0}{22.0}{6.0}
\drawcenteredtext{4.0}{22.0}{$A$}
\drawcenteredtext{50.0}{22.0}{$B$}
\end{picture}
}

\subsubsection{Parallel feedback}

The sequential feedback $\mathsf{Pfb}_{C}(\mathbf{Q}_{Y;A\times C,B\times
C}^{X})$ is represented by the diagram:

%
%
%
%
%
%
%
%
%
%
%
%
\centerline{\tt\setlength{\unitlength}{3.2pt}
\begin{picture}(56,48)
\thinlines
\drawframebox{27.0}{27.0}{22.0}{14.0}{$\mathbf{Q}$}
\drawdotline{26.0}{36.0}{26.0}{34.0}
\drawdotline{26.0}{20.0}{26.0}{18.0}
\drawpath{16.0}{30.0}{6.0}{30.0}
\drawpath{16.0}{24.0}{14.0}{24.0}
\drawpath{14.0}{24.0}{12.0}{22.0}
\drawpath{12.0}{22.0}{12.0}{16.0}
\drawpath{12.0}{16.0}{14.0}{14.0}
\drawpath{14.0}{14.0}{40.0}{14.0}
\drawpath{40.0}{14.0}{42.0}{16.0}
\drawpath{42.0}{16.0}{42.0}{22.0}
\drawpath{42.0}{22.0}{40.0}{24.0}
\drawpath{40.0}{24.0}{38.0}{24.0}
\drawpath{38.0}{30.0}{48.0}{30.0}
\drawframebox{27.0}{24.0}{34.0}{28.0}{}
\drawdotline{26.0}{42.0}{26.0}{38.0}
\drawdotline{26.0}{10.0}{26.0}{6.0}
\drawcenteredtext{4.0}{30.0}{$A$}
\drawcenteredtext{50.0}{30.0}{$B$}
\drawcenteredtext{14.0}{26.0}{$C$}
\drawcenteredtext{40.0}{26.0}{$C$}
\drawcenteredtext{28.0}{42.0}{$X$}
\drawcenteredtext{28.0}{6.0}{$Y$}
\end{picture}
}

\section{Examples}

\subsection{Any automaton is in $\Sigma(E)$}

If $E$ is the set of automata with two states with a single transition, then
any automaton may be given as a sequential expression of elements in $E$. We
illustrate by considering the first example of section 2. Let $T_{1}%
,T_{2},T_{3}$ be the three single transition automata as follows: $T_{1}$ has
the single transition labelled on the left by $a$ and the right by $(b_{1},c)$
with weight $2$, $T_{2}$ has the single transition labelled on the left by $a$
and the right by $(b_{1},c)$ with weight $3$. and $T_{3}$ has the single
transition labelled on the left by $a$ and the right by $(b_{2},c)$ with
weight $1$, \ Then the following diagram shows how the automaton may be given
as $T_{1}\mathbf{\boxplus} T_{2}\mathbf{\boxplus} T_{3}$ composed sequentially
with sequential wires:

%
%
%
%
%
%
%
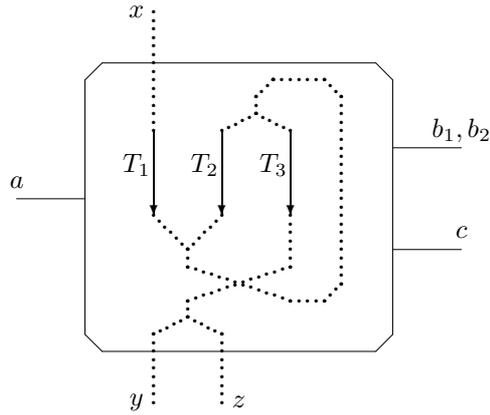
\centerline{\tt\setlength{\unitlength}{3.2pt}
\begin{picture}(66,58)
\thinlines
\drawpath{12.0}{44.0}{14.0}{46.0}
\drawpath{14.0}{46.0}{46.0}{46.0}
\drawpath{46.0}{46.0}{48.0}{44.0}
\drawpath{48.0}{44.0}{48.0}{14.0}
\drawpath{48.0}{14.0}{46.0}{12.0}
\drawpath{46.0}{12.0}{14.0}{12.0}
\drawpath{14.0}{12.0}{12.0}{14.0}
\drawpath{12.0}{14.0}{12.0}{44.0}
\drawpath{12.0}{30.0}{4.0}{30.0}
\drawpath{48.0}{36.0}{56.0}{36.0}
\drawpath{48.0}{24.0}{56.0}{24.0}
\drawcenteredtext{4.0}{32.0}{$a$}
\drawcenteredtext{56.0}{38.0}{$b_1,b_2$}
\drawcenteredtext{56.0}{26.0}{$c$}
\drawvector{20.0}{38.0}{10.0}{0}{-1}
\drawvector{28.0}{38.0}{10.0}{0}{-1}
\drawvector{36.0}{38.0}{10.0}{0}{-1}
\drawdotline{20.0}{38.0}{20.0}{52.0}
\drawdotline{32.0}{40.0}{28.0}{38.0}
\drawdotline{32.0}{40.0}{36.0}{38.0}
\drawdotline{32.0}{42.0}{32.0}{40.0}
\drawdotline{32.0}{42.0}{34.0}{44.0}
\drawdotline{34.0}{44.0}{34.0}{44.0}
\drawdotline{20.0}{28.0}{24.0}{24.0}
\drawdotline{24.0}{24.0}{28.0}{28.0}
\drawdotline{24.0}{24.0}{24.0}{22.0}
\drawdotline{36.0}{28.0}{36.0}{22.0}
\drawdotline{24.0}{22.0}{36.0}{18.0}
\drawdotline{36.0}{22.0}{24.0}{18.0}
\drawdotline{24.0}{18.0}{24.0}{16.0}
\drawdotline{34.0}{44.0}{40.0}{44.0}
\drawdotline{40.0}{44.0}{42.0}{42.0}
\drawdotline{42.0}{42.0}{42.0}{20.0}
\drawdotline{42.0}{20.0}{40.0}{18.0}
\drawdotline{40.0}{18.0}{36.0}{18.0}
\drawcenteredtext{18.0}{34.0}{$T_1$}
\drawcenteredtext{26.0}{34.0}{$T_2$}
\drawcenteredtext{34.0}{34.0}{$T_3$}
\drawdotline{24.0}{16.0}{20.0}{14.0}
\drawdotline{24.0}{16.0}{28.0}{14.0}
\drawdotline{28.0}{14.0}{28.0}{6.0}
\drawdotline{20.0}{14.0}{20.0}{6.0}
\drawcenteredtext{18.0}{52.0}{$x$}
\drawcenteredtext{18.0}{6.0}{$y$}
\drawcenteredtext{30.0}{6.0}{$z$}
\end{picture}
}

\subsection{The dining philosophers system}

The model of the dining philosophers problem we consider is an expression in
the algebra, involving also the automata $\mathbf{Phil}$ and $\mathbf{Fork}%
$.\ The system of $n$ dining philosophers is
\[
\mathbf{DF}_{n}=\mathsf{Pfb}\mathbf{(Phil}||\mathbf{Fork}|| \mathbf{Phil}%
||\mathbf{ Fork}||\cdots||\mathbf{Phil}||\mathbf{Fork}),
\]
where in this expression there are $n$ philosophers and $n$ forks.

The system is represented by the following diagram, where we abbreviate
$\mathbf{Phil}$ to $\mathcal{P}$ and $\mathbf{Fork}$ to $\mathcal{F}$.

%
%
%
%
%
%
%
%
%
%
%
%
%
%
%
%
%
%
%
%
%
%
%
%
%
%
\centerline{\tt\setlength{\unitlength}{4.5pt}
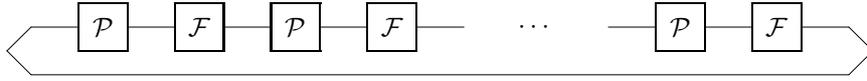
\begin{picture}(80,14)
\thinlines
\drawframebox{12.0}{8.0}{4.0}{4.0}{$\mathcal{P}$}
\drawframebox{20.0}{8.0}{4.0}{4.0}{$\mathcal{F}$}
\drawframebox{28.0}{8.0}{4.0}{4.0}{$\mathcal{P}$}
\drawframebox{36.0}{8.0}{4.0}{4.0}{$\mathcal{F}$}
\drawframebox{60.0}{8.0}{4.0}{4.0}{$\mathcal{P}$}
\drawframebox{68.0}{8.0}{4.0}{4.0}{$\mathcal{F}$}
\drawpath{14.0}{8.0}{18.0}{8.0}
\drawpath{22.0}{8.0}{26.0}{8.0}
\drawpath{30.0}{8.0}{34.0}{8.0}
\drawpath{38.0}{8.0}{42.0}{8.0}
\drawpath{54.0}{8.0}{58.0}{8.0}
\drawpath{62.0}{8.0}{66.0}{8.0}
\drawpath{70.0}{8.0}{74.0}{8.0}
\drawpath{74.0}{8.0}{76.0}{6.0}
\drawpath{76.0}{6.0}{74.0}{4.0}
\drawpath{74.0}{4.0}{6.0}{4.0}
\drawpath{6.0}{4.0}{4.0}{6.0}
\drawpath{4.0}{6.0}{6.0}{8.0}
\drawpath{6.0}{8.0}{10.0}{8.0}
\drawcenteredtext{48.0}{8.0}{$\cdots$}
\end{picture}}

Let us examine the case when $n=2$ with initial state $(1,1,1,1)$. Let
$\mathbf{Q}$ be the reachable part of $\mathbf{DF}_{2}$. The states reachable
from the initial state are $q_{1}=(1,1,1,1)$, $q_{2}=(1,3,3,2)$,
$q_{3}=(3,2,1,3)$, $q_{4}=(1,1,4,2)$, $q_{5}=(4,2,1,1)$, $q_{6}=(1,3,2,1)$,
$q_{7}=(2,1,1,3)$, $q_{8}=(2,3,2,3)$ ($q_{8}$ is the unique deadlock state).
The single matrix of the automaton $\mathsf{Q}$, using this ordering of the
states, is
\[
\left[
\begin{array}
[c]{cccccccc}%
\frac{1}{4} & 0 & 0 & 0 & 0 & \frac{1}{4} & \frac{1}{4} & \frac{1}{4}\\
0 & \frac{1}{2} & 0 & \frac{1}{2} & 0 & 0 & 0 & 0\\
0 & 0 & \frac{1}{2} & 0 & \frac{1}{2} & 0 & 0 & 0\\
\frac{1}{2} & 0 & 0 & \frac{1}{2} & 0 & 0 & 0 & 0\\
\frac{1}{2} & 0 & 0 & 0 & \frac{1}{2} & 0 & 0 & 0\\
0 & \frac{1}{3} & 0 & 0 & 0 & \frac{1}{3} & 0 & \frac{1}{3}\\
0 & 0 & \frac{1}{3} & 0 & 0 & 0 & \frac{1}{3} & \frac{1}{3}\\
0 & 0 & 0 & 0 & 0 & 0 & 0 & 1
\end{array}
\right]  .
\]

Calculating powers of this matrix we see that the probability of reaching
deadlock from the initial state in $2$ steps is $\frac{23}{48}$, in $3$ steps
is $\frac{341}{576}$, and in $4$ steps is $\frac{4415}{6912}$.

\subsubsection{ The probability of deadlock}

We describe here the result of \cite{dFSW09a} in which we apply
Perron-Frobenius theory to the problem of deadlock in the Dining Philosopher problem..

\begin{dfn}
Consider a Markov automaton $\mathbf{Q}$ with input and output sets being one
element sets $\{\varepsilon\}$. A state $q$ is called a \emph{deadlock} if the
only transition out of $q$ with positive probability is a transition from $q$
to $q$ (the probability of the transition must necessarily be $1$).
\end{dfn}

\begin{prop}
(Perron-Frobenius) Consider a Markov automaton $\mathbf{Q}$ with interfaces
being one element sets, with an initial state $q_{0}$. Suppose that (i)
$\mathbf{Q}$ has precisely one reachable deadlock state, (ii) for each
reachable state, not a deadlock, there is a path with non-zero probability to
$q_{0}$, and (iii) for each reachable state $q$ there is a transition with
non-zero probability to itself.

\noindent Then the probability of reaching a deadlock from the initial state
in $k$ steps tends to $1$ as $k$ tends to infinity.
\end{prop}

\begin{cor}
In the dining philosopher problem $\mathbf{DF}_{n}$ with $q_{0}$ being the
state $(1,1,\cdots,1)$ the unique reachable deadlock is $(3,2,3,2,\cdots
,3,2)$. The initial state is reachable from all other reachable states. Hence
the probability of reaching a deadlock from the initial state in $k$ steps
tends to $1$ as $k$ tends to infinity.
\end{cor}

\noindent\textbf{Remark. }The corollary does not depend on the specific
positive probabilities of the actions of the philosophers and forks. Hence the
result is true with any positive probabilities. In fact, different
philosophers and forks may have different probabilities without affecting the
conclusion of the corollary.

\section{Sofia's birthday party}

The example we would like to describe is a variant of the Dining Philosopher
Problem which we call Sofia's Birthday Party. Instead of a circle of
philosophers around a table with as many forks, we consider a circle of seats
around a table separated by forks on the table. Then there are a number of
children (not greater than the number of seats). The protocol of each child is
the same as that of a philosopher. However in addition, if a child is not
holding a fork, and the seat to the right is empty, the child may change seats
- the food may be better there.

\emph{To simplify the problem we will assume that all transitions have weight
}$0$ \emph{or} $1$, \emph{so the transitions of components we mention will all
have weight}$1$.

To describe this we need six automata -- a child $\mathcal{C}$, an empty seat
$\mathcal{E}$, a fork $\mathcal{F}$, two transition elements $\mathcal{L}$ and
$\mathcal{R}$, and the identity $1_{A}$ of $A$ (a wire). The interface sets
involved are $A=\{x,\varepsilon\}$ and $B=\{\varepsilon,t,r\}$.

The transition elements have left and right interfaces $A\times B.$ The graph
of the of the transition element $\mathcal{L}$ has two vertices $p$ and $q$
and one labelled edges $x,\varepsilon/\varepsilon,\varepsilon:q\rightarrow p.$
Its top interface is $Q=\{q\}$, and its bottom interface is $P=\{p\}$. The
graph of the transition element $\mathcal{R}$ also has two vertices $p$ and
$q$, and has one labelled edges $\varepsilon,\varepsilon/x,\varepsilon
:q\rightarrow p$. Its top interface is also $Q=\{q\}$, and its bottom
interface is $P=\{p\}$. The empty seat $\mathcal{E}$ has left and right
interfaces $A\times B$. The graph  of the empty seat has one vertex $e$ and
one labelled edge $\varepsilon,\varepsilon/\varepsilon,\varepsilon
:e\rightarrow e$. Its top interface is $P$ and its bottom interface is $Q$.
The functions $\gamma_{0},\gamma_{1}$ are uniquely defined.

The child $\mathcal{C}$ has labelled graph as follows:

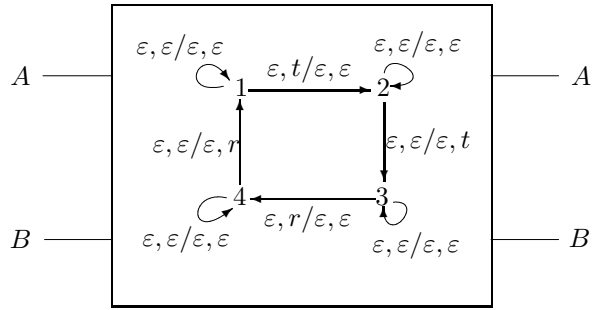
\centerline{
{\tt\setlength{\unitlength}{2.0pt}
\begin{picture}(115,65)
\thinlines
\drawcenteredtext{46.7}{46.16}{$1$}
\drawcenteredtext{73.29}{46.18}{$2$}
\drawcenteredtext{73.01}{25.69}{$3$}
\drawcenteredtext{46.43}{25.68}{$4$}
\drawvector{48.11}{45.61}{22.78}{1}{0}
\drawvector{73.43}{43.36}{14.86}{0}{-1}
\drawvector{71.75}{25.13}{23.34}{-1}{0}
\drawvector{46.43}{27.93}{15.98}{0}{1}
\path
(43.34,46.44)(43.34,46.44)(43.06,46.38)(42.79,46.33)(42.54,46.3)(42.29,46.26)(42.04,46.23)(41.81,46.23)(41.59,46.2)(41.38,46.2)
\path
(41.38,46.2)(41.16,46.2)(40.97,46.22)(40.77,46.23)(40.59,46.25)(40.41,46.26)(40.25,46.3)(40.09,46.33)(39.93,46.37)(39.79,46.41)
\path
(39.79,46.41)(39.65,46.47)(39.52,46.52)(39.4,46.58)(39.27,46.65)(39.16,46.7)(39.06,46.77)(38.97,46.86)(38.88,46.93)(38.79,47.01)
\path
(38.79,47.01)(38.72,47.08)(38.65,47.18)(38.59,47.26)(38.52,47.34)(38.47,47.44)(38.43,47.54)(38.38,47.62)(38.34,47.73)(38.33,47.83)
\path
(38.33,47.83)(38.29,47.93)(38.27,48.02)(38.27,48.12)(38.25,48.23)(38.25,48.33)(38.27,48.43)(38.27,48.52)(38.29,48.62)(38.31,48.73)
\path
(38.31,48.73)(38.33,48.83)(38.36,48.93)(38.38,49.02)(38.43,49.12)(38.47,49.22)(38.52,49.3)(38.56,49.41)(38.61,49.5)(38.68,49.58)
\path
(38.68,49.58)(38.74,49.66)(38.81,49.75)(38.88,49.83)(38.95,49.91)(39.04,49.98)(39.11,50.05)(39.2,50.12)(39.29,50.18)(39.38,50.23)
\path
(39.38,50.23)(39.47,50.29)(39.58,50.33)(39.68,50.38)(39.77,50.43)(39.88,50.45)(40.0,50.48)(40.11,50.51)(40.22,50.54)(40.34,50.55)
\path
(40.34,50.55)(40.47,50.55)(40.59,50.55)(40.7,50.55)(40.84,50.54)(40.97,50.51)(41.09,50.5)(41.22,50.47)(41.36,50.43)(41.5,50.37)
\path
(41.5,50.37)(41.63,50.33)(41.77,50.26)(41.91,50.19)(42.06,50.11)(42.2,50.02)(42.34,49.93)(42.5,49.83)(42.65,49.7)(42.79,49.58)
\path
(42.79,49.58)(42.93,49.45)(43.09,49.3)(43.24,49.16)(43.38,49.0)(43.54,48.83)(43.7,48.63)(43.84,48.44)(44.0,48.25)(44.15,48.02)
\path(44.15,48.02)(44.29,47.8)(44.45,47.56)(44.45,47.56)
\drawvector{43.7}{48.63}{0.75}{3}{-4}
\path
(73.43,48.97)(73.43,48.97)(73.44,49.11)(73.48,49.25)(73.52,49.38)(73.56,49.51)(73.62,49.62)(73.66,49.73)(73.73,49.84)(73.79,49.94)
\path
(73.79,49.94)(73.86,50.04)(73.93,50.12)(74.0,50.2)(74.08,50.27)(74.16,50.34)(74.25,50.41)(74.33,50.47)(74.43,50.51)(74.51,50.56)
\path
(74.51,50.56)(74.61,50.61)(74.7,50.63)(74.8,50.66)(74.91,50.69)(75.01,50.7)(75.12,50.73)(75.23,50.73)(75.33,50.73)(75.44,50.73)
\path
(75.44,50.73)(75.56,50.73)(75.68,50.73)(75.79,50.72)(75.9,50.69)(76.01,50.68)(76.12,50.66)(76.23,50.62)(76.36,50.59)(76.47,50.55)
\path
(76.47,50.55)(76.58,50.51)(76.69,50.48)(76.8,50.44)(76.91,50.38)(77.01,50.33)(77.12,50.29)(77.23,50.23)(77.33,50.16)(77.44,50.11)
\path
(77.44,50.11)(77.54,50.05)(77.63,49.98)(77.73,49.91)(77.81,49.84)(77.91,49.77)(78.0,49.69)(78.08,49.62)(78.16,49.55)(78.23,49.48)
\path
(78.23,49.48)(78.3,49.4)(78.37,49.31)(78.44,49.23)(78.5,49.16)(78.55,49.06)(78.61,48.98)(78.65,48.91)(78.69,48.81)(78.73,48.73)
\path
(78.73,48.73)(78.76,48.65)(78.79,48.56)(78.8,48.48)(78.81,48.4)(78.83,48.31)(78.83,48.23)(78.81,48.15)(78.8,48.06)(78.79,47.98)
\path
(78.79,47.98)(78.76,47.9)(78.73,47.81)(78.69,47.73)(78.65,47.66)(78.58,47.58)(78.52,47.51)(78.45,47.44)(78.37,47.37)(78.3,47.3)
\path
(78.3,47.3)(78.19,47.23)(78.09,47.16)(77.98,47.09)(77.86,47.04)(77.73,46.98)(77.58,46.91)(77.44,46.87)(77.29,46.81)(77.12,46.76)
\path
(77.12,46.76)(76.94,46.72)(76.75,46.68)(76.55,46.63)(76.33,46.59)(76.12,46.56)(75.87,46.54)(75.63,46.51)(75.38,46.48)(75.12,46.47)
\path(75.12,46.47)(74.83,46.45)(74.55,46.44)(74.55,46.44)
\drawvector{75.63}{46.51}{1.08}{-1}{0}
\path
(75.12,24.56)(75.12,24.56)(75.3,24.58)(75.48,24.59)(75.66,24.59)(75.81,24.58)(75.98,24.56)(76.12,24.56)(76.26,24.54)(76.4,24.5)
\path
(76.4,24.5)(76.52,24.47)(76.65,24.45)(76.76,24.4)(76.86,24.36)(76.95,24.31)(77.05,24.27)(77.12,24.2)(77.2,24.15)(77.27,24.09)
\path
(77.27,24.09)(77.33,24.02)(77.4,23.95)(77.44,23.88)(77.48,23.81)(77.52,23.72)(77.55,23.65)(77.58,23.56)(77.61,23.49)(77.62,23.4)
\path
(77.62,23.4)(77.63,23.31)(77.65,23.22)(77.65,23.13)(77.65,23.04)(77.63,22.95)(77.62,22.84)(77.61,22.75)(77.58,22.65)(77.55,22.56)
\path
(77.55,22.56)(77.52,22.47)(77.5,22.36)(77.45,22.27)(77.41,22.18)(77.37,22.09)(77.33,21.99)(77.27,21.9)(77.22,21.79)(77.16,21.7)
\path
(77.16,21.7)(77.09,21.61)(77.04,21.52)(76.97,21.45)(76.91,21.36)(76.83,21.27)(76.76,21.2)(76.69,21.11)(76.61,21.04)(76.52,20.97)
\path
(76.52,20.97)(76.44,20.9)(76.37,20.83)(76.29,20.75)(76.19,20.7)(76.11,20.65)(76.02,20.59)(75.94,20.54)(75.84,20.5)(75.76,20.45)
\path
(75.76,20.45)(75.66,20.41)(75.58,20.38)(75.48,20.36)(75.4,20.33)(75.3,20.31)(75.22,20.29)(75.12,20.29)(75.04,20.29)(74.94,20.29)
\path
(74.94,20.29)(74.86,20.31)(74.76,20.33)(74.68,20.34)(74.58,20.38)(74.51,20.4)(74.43,20.45)(74.33,20.5)(74.26,20.56)(74.18,20.63)
\path
(74.18,20.63)(74.09,20.7)(74.02,20.77)(73.95,20.86)(73.87,20.95)(73.8,21.06)(73.75,21.18)(73.68,21.29)(73.62,21.43)(73.55,21.56)
\path
(73.55,21.56)(73.51,21.7)(73.45,21.86)(73.41,22.02)(73.36,22.2)(73.31,22.38)(73.27,22.59)(73.23,22.79)(73.2,23.0)(73.19,23.22)
\path(73.19,23.22)(73.16,23.47)(73.15,23.7)(73.15,23.72)
\drawvector{73.27}{22.59}{1.13}{0}{1}
\path
(43.9,25.4)(43.9,25.4)(43.68,25.45)(43.47,25.49)(43.25,25.52)(43.04,25.54)(42.84,25.56)(42.65,25.58)(42.47,25.58)(42.27,25.58)
\path
(42.27,25.58)(42.09,25.56)(41.91,25.56)(41.75,25.54)(41.58,25.52)(41.41,25.49)(41.25,25.45)(41.09,25.4)(40.95,25.36)(40.81,25.31)
\path
(40.81,25.31)(40.65,25.27)(40.52,25.2)(40.4,25.15)(40.27,25.08)(40.13,25.0)(40.02,24.93)(39.9,24.86)(39.79,24.77)(39.7,24.7)
\path
(39.7,24.7)(39.59,24.61)(39.5,24.52)(39.4,24.43)(39.31,24.34)(39.24,24.25)(39.15,24.15)(39.09,24.04)(39.02,23.95)(38.95,23.84)
\path
(38.95,23.84)(38.9,23.75)(38.84,23.65)(38.79,23.54)(38.74,23.43)(38.7,23.34)(38.66,23.22)(38.63,23.13)(38.61,23.02)(38.59,22.91)
\path
(38.59,22.91)(38.56,22.81)(38.56,22.72)(38.54,22.61)(38.54,22.5)(38.54,22.4)(38.54,22.31)(38.56,22.22)(38.58,22.13)(38.59,22.02)
\path
(38.59,22.02)(38.63,21.93)(38.65,21.86)(38.68,21.77)(38.72,21.68)(38.77,21.61)(38.81,21.54)(38.88,21.45)(38.93,21.38)(39.0,21.33)
\path
(39.0,21.33)(39.06,21.27)(39.13,21.2)(39.22,21.15)(39.29,21.11)(39.38,21.06)(39.47,21.02)(39.58,20.99)(39.68,20.95)(39.77,20.93)
\path
(39.77,20.93)(39.9,20.91)(40.0,20.9)(40.13,20.9)(40.25,20.9)(40.38,20.9)(40.52,20.91)(40.66,20.93)(40.81,20.97)(40.95,21.0)
\path
(40.95,21.0)(41.11,21.04)(41.27,21.09)(41.45,21.15)(41.61,21.2)(41.79,21.27)(41.97,21.36)(42.15,21.45)(42.34,21.54)(42.54,21.65)
\path
(42.54,21.65)(42.75,21.75)(42.95,21.88)(43.15,22.02)(43.38,22.15)(43.59,22.31)(43.81,22.47)(44.04,22.65)(44.29,22.83)(44.52,23.02)
\path(44.52,23.02)(44.77,23.22)(45.02,23.43)(45.02,23.43)
\drawvector{43.81}{22.47}{1.2}{4}{3}
\drawcenteredtext{59.37}{49.52}{$\varepsilon,t/\varepsilon,\varepsilon$}
\drawcenteredtext{81.58}{36.06}{$\varepsilon, \varepsilon/\varepsilon, t$}
\drawcenteredtext{59.08}{21.2}{$\varepsilon,r/\varepsilon,\varepsilon$}
\drawcenteredtext{38.27}{35.22}{$\varepsilon,\varepsilon/\varepsilon,r$}
\drawcenteredtext{35.2}{53.73}{$\varepsilon,\varepsilon/\varepsilon
,\varepsilon$}
\drawcenteredtext{79.64}{54.3}{$\varepsilon,\varepsilon/\varepsilon
,\varepsilon$}
\drawcenteredtext{79.36}{16.43}{$\varepsilon,\varepsilon/\varepsilon
,\varepsilon$}
\drawcenteredtext{36.32}{17.54}{$\varepsilon,\varepsilon/\varepsilon
,\varepsilon$}
\drawframebox{57.96}{33.26}{70.87}{56.66}{}
\drawpath{22.53}{48.41}{9.59}{48.41}
\drawpath{22.53}{17.54}{9.87}{17.54}
\drawpath{93.4}{48.41}{106.05}{48.41}
\drawpath{93.4}{17.54}{106.05}{17.54}
\drawcenteredtext{5.36}{48.41}{$A$}
\drawcenteredtext{110.27}{48.41}{$A$}
\drawcenteredtext{110.0}{17.84}{$B$}
\drawcenteredtext{5.36}{17.55}{$B$}
\end{picture}
}}

The states have the following interpretation: in state $1$ the child has no
forks; in state $2$ it has a fork in its left hand; in state $3$ it has both
forks (and can eat); in state $4$ it has returned it left fork. The child's
top interface is $P$ and its bottom interface is $Q$. The function $\gamma
_{0}$ takes $p$ to $1;$ the function $\gamma_{1}$ takes $q$ to $1$.

The fork $\mathcal{F}$ is as in the dining philosopher system (but with all
transtions weighted $1$).

Let $\mathcal{S=}\mathsf{Cfb}_{P}(C\bullet R\bullet E\bullet L)$. This
automaton has the following interpretation -- it can either be a child (on a
seat) or an empty seat. The transition elements $\mathcal{R}$ and
$\mathcal{L}$ allow the seat to become occupied or vacated. It is
straightforward to see that this automaton is a positive weighted automaton.

Then Sofia's Birthday Party is given by the normalization of expression%
\[
\mathsf{Pfb}_{A\times B}(\mathcal{S||(}1_{A}\times\mathcal{F)||S||(}%
1_{A}\times\mathcal{F)||...||S||(}1_{X}\times\mathcal{F)).}%
\]

This automaton has the behaviour as informally described above. Its
diagrammatic representation is:

%
%
%
%
%
%
%
%
%
%
%
%
%
%
%
%
%
%
%
%
%
%
%
%
%
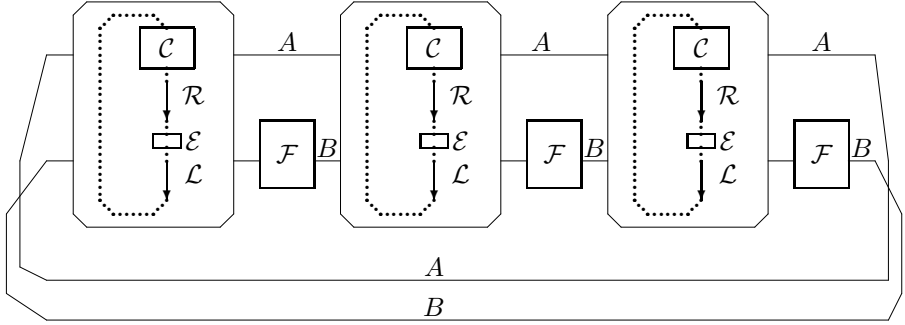
\centerline{
{\tt\setlength{\unitlength}{2.5pt}
\begin{picture}(142,56)
\thinlines
\drawpath{42.0}{28.0}{38.0}{28.0}
\drawpath{16.0}{52.0}{36.0}{52.0}
\drawframebox{28.0}{45.0}{8.0}{6.0}{$\mathcal{C}$}
\drawvector{28.0}{40.0}{6.0}{0}{-1}
\drawframebox{28.0}{31.0}{4.0}{2.0}{}
\drawcenteredtext{32}{31}{$\mathcal{E}$}
\drawvector{28.0}{28.0}{6.0}{0}{-1}
\drawdotline{28.0}{48.0}{26.0}{50.0}
\drawdotline{26.0}{50.0}{20.0}{50.0}
\drawdotline{20.0}{50.0}{18.0}{48.0}
\drawdotline{28.0}{42.0}{28.0}{42.0}
\drawdotline{28.0}{42.0}{28.0}{40.0}
\drawdotline{28.0}{34.0}{28.0}{32.0}
\drawdotline{28.0}{30.0}{28.0}{28.0}
\drawdotline{18.0}{48.0}{18.0}{22.0}
\drawdotline{18.0}{22.0}{20.0}{20.0}
\drawdotline{28.0}{22.0}{26.0}{20.0}
\drawdotline{26.0}{20.0}{20.0}{20.0}
\drawpath{36.0}{52.0}{38.0}{50.0}
\drawpath{38.0}{50.0}{38.0}{20.0}
\drawpath{38.0}{20.0}{36.0}{18.0}
\drawpath{36.0}{18.0}{16.0}{18.0}
\drawpath{16.0}{18.0}{14.0}{20.0}
\drawpath{14.0}{20.0}{14.0}{50.0}
\drawpath{14.0}{50.0}{16.0}{52.0}
\drawcenteredtext{32.0}{38.0}{$\mathcal{R}$}
\drawframebox{46.0}{29.0}{8.0}{10.0}{$\mathcal{F}$}
\drawframebox{86.0}{29.0}{8.0}{10.0}{$\mathcal{F}$}
\drawframebox{126.0}{29.0}{8.0}{10.0}{$\mathcal{F}$}
\drawpath{38.0}{44.0}{54.0}{44.0}
\drawpath{56.0}{52.0}{76.0}{52.0}
\drawpath{50.0}{28.0}{54.0}{28.0}
\drawpath{78.0}{28.0}{82.0}{28.0}
\drawpath{90.0}{28.0}{94.0}{28.0}
\drawpath{118.0}{28.0}{122.0}{28.0}
\drawpath{130.0}{28.0}{134.0}{28.0}
\drawpath{14.0}{28.0}{10.0}{28.0}
\drawpath{14.0}{44.0}{10.0}{44.0}
\drawpath{78.0}{44.0}{94.0}{44.0}
\drawpath{118.0}{44.0}{134.0}{44.0}
\drawpath{10.0}{10.0}{136.0}{10.0}
\drawpath{10.0}{4.0}{136.0}{4.0}
\drawpath{10.0}{28.0}{4.0}{20.0}
\drawpath{4.0}{20.0}{4.0}{8.0}
\drawpath{4.0}{8.0}{10.0}{4.0}
\drawpath{10.0}{44.0}{6.0}{28.0}
\drawpath{6.0}{28.0}{6.0}{12.0}
\drawpath{6.0}{12.0}{10.0}{10.0}
\drawpath{134.0}{28.0}{138.0}{20.0}
\drawpath{138.0}{20.0}{138.0}{8.0}
\drawpath{138.0}{8.0}{136.0}{4.0}
\drawpath{134.0}{44.0}{136.0}{28.0}
\drawpath{136.0}{28.0}{136.0}{10.0}
\drawcenteredtext{46.0}{46.0}{$A$}
\drawcenteredtext{52.0}{30.0}{$B$}
\drawcenteredtext{84.0}{46.0}{$A$}
\drawcenteredtext{92.0}{30.0}{$B$}
\drawcenteredtext{126.0}{46.0}{$A$}
\drawcenteredtext{132.0}{30.0}{$B$}
\drawcenteredtext{68.0}{12.0}{$A$}
\drawcenteredtext{68.0}{6.0}{$B$}
\drawframebox{68.0}{45.0}{8.0}{6.0}{$\mathcal{C}$}
\drawvector{68.0}{40.0}{6.0}{0}{-1}
\drawframebox{68.0}{31.0}{4.0}{2.0}{}
\drawcenteredtext{72}{31}{$\mathcal{E}$}
\drawvector{68.0}{28.0}{6.0}{0}{-1}
\drawdotline{68.0}{48.0}{66.0}{50.0}
\drawdotline{66.0}{50.0}{60.0}{50.0}
\drawdotline{60.0}{50.0}{58.0}{48.0}
\drawdotline{68.0}{42.0}{68.0}{42.0}
\drawdotline{68.0}{42.0}{68.0}{40.0}
\drawdotline{68.0}{34.0}{68.0}{32.0}
\drawdotline{68.0}{30.0}{68.0}{28.0}
\drawdotline{58.0}{48.0}{58.0}{22.0}
\drawdotline{58.0}{22.0}{60.0}{20.0}
\drawdotline{68.0}{22.0}{66.0}{20.0}
\drawdotline{66.0}{20.0}{60.0}{20.0}
\drawpath{76.0}{52.0}{78.0}{50.0}
\drawpath{78.0}{50.0}{78.0}{20.0}
\drawpath{78.0}{20.0}{76.0}{18.0}
\drawpath{76.0}{18.0}{56.0}{18.0}
\drawpath{56.0}{18.0}{54.0}{20.0}
\drawpath{54.0}{20.0}{54.0}{50.0}
\drawpath{54.0}{50.0}{56.0}{52.0}
\drawcenteredtext{72.0}{38.0}{$\mathcal{R}$}
\drawpath{96.0}{52.0}{116.0}{52.0}
\drawframebox{108.0}{45.0}{8.0}{6.0}{$\mathcal{C}$}
\drawvector{108.0}{40.0}{6.0}{0}{-1}
\drawframebox{108.0}{31.0}{4.0}{2.0}{}
\drawcenteredtext{112}{31}{$\mathcal{E}$}\drawvector{108.0}{28.0}{6.0}{0}{-1}
\drawdotline{108.0}{48.0}{106.0}{50.0}
\drawdotline{106.0}{50.0}{100.0}{50.0}
\drawdotline{100.0}{50.0}{98.0}{48.0}
\drawdotline{108.0}{42.0}{108.0}{42.0}
\drawdotline{108.0}{42.0}{108.0}{40.0}
\drawdotline{108.0}{34.0}{108.0}{32.0}
\drawdotline{108.0}{30.0}{108.0}{28.0}
\drawdotline{98.0}{48.0}{98.0}{22.0}
\drawdotline{98.0}{22.0}{100.0}{20.0}
\drawdotline{108.0}{22.0}{106.0}{20.0}
\drawdotline{106.0}{20.0}{100.0}{20.0}
\drawpath{116.0}{52.0}{118.0}{50.0}
\drawpath{118.0}{50.0}{118.0}{20.0}
\drawpath{118.0}{20.0}{116.0}{18.0}
\drawpath{116.0}{18.0}{96.0}{18.0}
\drawpath{96.0}{18.0}{94.0}{20.0}
\drawpath{94.0}{20.0}{94.0}{50.0}
\drawpath{94.0}{50.0}{96.0}{52.0}
\drawcenteredtext{112.0}{38.0}{$\mathcal{R}$}
\drawcenteredtext{32.0}{26.0}{$\mathcal{L}$}
\drawcenteredtext{72.0}{26.0}{$\mathcal{L}$}
\drawcenteredtext{112.0}{26.0}{$\mathcal{L}$}
\end{picture}
}}

Notice that though Sofia's Birthday Party belongs to $\Pi\Sigma\Sigma(E)$
slight variations of the system belong instead to $\Pi\Sigma\Pi\Sigma(E)$, for
example the system where more than one child may occupy a seat (communicating
there). If the system starts in a state with as many children as seats - then
movement is impossible and the system is equivalent to the dining philosophers.

Let us look at a particular case of Sofia's Birthday Party in more detail.
Consider the system with three seats, and two children. There are $36$ states
reachable from initial state $(5,1,1,1,1,1)$, where $5$ is the state in which
the seat is empty, and there are $141$ transitions.

The states are%
\[%
\begin{array}
[c]{cccc}%
(5,1,1,1,1,1) & (5,1,1,3,2,1) & (5,3,2,1,1,1) & (5,3,2,3,2,1)\\
(1,1,1,1,5,1) & (1,3,2,1,5,1) & (5,1,1,3,3,2) & (5,3,2,3,3,2)\\
(5,3,3,2,1,1) & (1,3,3,2,5,1) & (1,1,5,1,1,1) & (2,1,1,1,5,3)\\
(2,1,5,1,1,3) & (2,3,2,1,5,3) & (2,3,3,2,5,3) & (5,1,1,1,4,2)\\
(5,3,2,1,4,2) & (5,1,4,2,1,1) & (1,1,4,2,5,1) & (2,1,4,2,5,3)\\
(1,1,5,3,2,1) & (2,1,5,3,2,3) & (3,2,1,1,5,3) & (3,2,5,1,1,3)\\
(3,2,5,3,2,3) & (5,3,3,2,4,2) & (3,2,4,2,5,3) & (1,1,5,3,3,2)\\
(4,2,1,1,5,1) & (4,2,5,1,1,1) & (4,2,5,3,2,1) & (5,1,4,2,4,2)\\
(4,2,4,2,5,1) & (1,1,5,1,4,2) & (4,2,5,3,3,2) & (4,2,5,1,4,2).
\end{array}
\]

Then in $12$ of these states there is a child eating: only one may eat at a
time. It is straightforward to calculate the $36\times36$ Markov matrix and
iterating show that the probability of a child eating after $1$ step from
initial state $(5,1,1,1,1,1)$ is $0$, after $2$ steps is $\frac{19}{60}$,
after $3$ steps is $\frac{98}{225}$, after $4$ steps is $\frac{49133}{108000}%
$, after $5$ steps is $\frac{1473023}{3240000}$ and after $100$ steps is
$0.3768058221$.

\end{document}